\documentclass[aos,preprint]{imsart}

\usepackage{fullpage}

\usepackage{latexsym,amsmath,amssymb,amsthm,epic,eepic,multirow}
\usepackage{natbib}

\usepackage{pgfplots}

\usepackage{graphicx}
\usepackage{multirow}
\usepackage{upgreek}
\usepackage{subfigure}
\usepackage{algorithm}
\usepackage[normalem]{ulem}
\usepackage{algorithmic}
\usepackage{enumitem}
\usepackage{comment}
\usepackage{mathrsfs}
\usepackage{dsfont}

\usepackage{float}
\usepackage{hyperref}
   \hypersetup{ colorlinks = true, linkcolor = blue, citecolor=blue }

\allowdisplaybreaks[4]

\newtheoremstyle{newstyle}
   {\topsep}{\topsep}{\itshape}{}{\bfseries}{.}{.5em}{{\thmname{#1 }}{\thmnumber{#2}}{\thmnote{ (#3)}}}
\theoremstyle{newstyle}
\newtheorem{theorem}{Theorem}
\newtheorem{proposition}{Proposition}
\newtheorem{lemma}{Lemma}
\newtheorem{corollary}{Corollary}
\newtheorem{condition}{Condition}

\theoremstyle{definition}
\newtheorem{example}{Example}

\theoremstyle{remark}
\newtheorem{remark}{Remark}

\setattribute{journal}{name}{}

\begin{document}

\begin{frontmatter}
\title{Slightly Conservative Bootstrap for Maxima of Sums}

\runtitle{Slightly Conservative Bootstrap}

\begin{aug}
\author{\fnms{Hang} \snm{Deng}\ead[label=e1]{hdeng@stat.rutgers.edu}}

\runauthor{Hang Deng}

\thankstext{t1}{Partially supported by DMS-1454817 and CCF-1934924.}

\affiliation{Department of Statistics, Rutgers University}

\end{aug}

\maketitle

\begin{abstract}
We study the bootstrap for the maxima of the sums of independent random variables, a problem of high relevance to many applications in modern statistics. Since the consistency of bootstrap was justified by Gaussian approximation in \cite{chernozhukov2013}, quite a few attempts 
have been made to sharpen the error bound for bootstrap and reduce the sample size requirement for bootstrap consistency. In this paper, we show that the sample size requirement can be dramatically improved when we make the inference slightly conservative, that is, to inflate the bootstrap quantile $t_{\alpha}^*$ by a small fraction, e.g. by $1\%$ to $1.01\,t^*_\alpha$. This simple procedure yields error bounds for the coverage probability of conservative bootstrap at as fast a rate as $\sqrt{(\log p)/n}$ under suitable conditions, so that not only the sample size requirement can be reduced to $\log p \ll n$ but also the overall convergence rate is nearly parametric. Furthermore, we improve the error bound for the coverage probability of the standard non-conservative bootstrap to $[(\log (np))^3 (\log p)^2/n]^{1/4}$ under general assumptions on data. These results are established for the empirical bootstrap and the multiplier bootstrap with third moment match. An improved coherent Lindeberg interpolation method, originally proposed in \cite{deng2017beyond}, is developed to derive sharper comparison bounds, especially for the maxima. 
\end{abstract}

\end{frontmatter}






\section{Introduction} Let $ \boldsymbol{X} =(X_1,\ldots, X_n)^\top \in \mathbb{R}^{n\times p}$ be a random matrix with
independent rows $X_i = (X_{i,1},\ldots,X_{i,p})^\top \in \mathbb{R}^p$, $1 \le i \le n$
. We are interested in the approximation of the distribution of the maximum of normalized sums
\begin{align}\label{T_n} 
T_{n} = \max_{1\le j \le p} \frac{1}{\sqrt{n}}\sum_{i=1}^n \big( X_{i,j} - \mathbb{E} \overline{X}_{n,j} \big) \hbox{ where } \overline{X}_{n,j} = \frac{1}{n}\sum_{i=1}^{n}X_{i,j}. 
\end{align}
We focus on the case of large $p= p_n$, including exponential growth of $p_n$ at certain rate as $n\to\infty$. 

The approximation of the distribution of $T_n$ has broad applications.  
Examples include sure screening \citep{fanlv07}, removing spurious correlation \citep{fan2016guarding},
testing the equality of two matrices \citep{cai2013two} and \citep{chang2016comparing}, 
detecting ridges and estimating level sets \citep{chen2015asymptotic,chen2016density}, among others.
It can be also used in time series settings \citep{zhang2017gaussian} and high-dimensional regression \citep{zhangzhang11,belloni2014inference,bellonietal15,ZhangCheng2016,dezeure2016high}.
In such modern applications, $p = p_n$ is not fixed and typically much larger than $n$. 

The problem can be put into a broader context as a subproblem of the approximation of $\mathbb{P} \{ \sum_{i=1}^n ( X_i - \mathbb{E} \overline{X}_n )/\sqrt{n} \in A \}$, where $A$ belongs to a certain collection of sets $\mathcal{A}$. \cite{bentkus2003dependence} studied the general Gaussian approximation in which $\mathcal{A}$ is the collection of all convex sets and gave a consistency result that essentially requires $p^{7/2} \ll n$, and \cite{zhilova2016non} improved upon the result 
on the set of Euclidean balls $\mathcal{A} := \big\{ \{x \in \mathbb{R}^p: \|x\|_2 \le t\}, t >0 \big\}$. 

\cite{chernozhukov2013} first studied the approximation of the distribution function of $T_n$ in Kolmogorov-Smirnov distance and used the Gaussian approximation to establish the consistency of the empirical bootstrap and the Gaussian multiplier bootstrap for $T_n$, 
\begin{align}\label{intro:rate}
\Big| (1- \alpha) - \mathbb{P} \big\{ T_n \le t_{\alpha}^* \big\} \Big| \le C \Big(\frac{(\log (np))^7}{n}\Big)^{1/8},
\end{align}
under certain moment and tail probability conditions on the data $\{X_i\}$. Here $C$ is a fixed constant and $t_{\alpha}^*$ is the $(1- \alpha)$-quantile of the bootstrap version of $T_n$. Later in \cite{chernozhukov2017central} this result was sharpened to $C\big((\log (np))^7/n \big)^{1/6}$; general hyperrectangular sets and sparse convex sets $\mathcal{A}$ were also studied therein. \cite{deng2017beyond} went beyond Gaussian approximation and used a coherent Lindeberg interpolation to improve the rate in \eqref{intro:rate} to $C\big((\log (np))^5/n \big)^{1/6}$ for the empirical bootstrap or the multiplier bootstrap with third moment match \citep{liu1988bootstrap,mammenwild1993}. 
In the sequel, the use of the bootstrapped $t_{\alpha}^*$ to estimate the $(1- \alpha)-$quantile of $T_n$ is referred to as the \textit{exact bootstrap}. For such methods, we aim to control the two-sided error as in \eqref{intro:rate}. 

Although the error bound has been improved to $C\big((\log (np))^5/n \big)^{1/6}$, we are still naturally interested in further weakening the sample size requirement $n\gg (\log p)^5$ and improving the rate of consistency towards the parametric rate $n^{-1/2}$. Our main contributions of this paper are two-fold: (i) We study a slightly \textit{conservative bootstrap} procedure for which we substantially improve the convergence rate to as fast as $(( \log p)/n )^{1/2}$ and therefore the sample size requirement to $n \gg \log p$; (ii) We improve the convergence rate for the exact bootstrap to $\big((\log(np))^3(\log p)^2/n \big)^{1/4}$. Moreover, we establish these results under general moment and tail probability conditions on the data $\{X_i\}$ so that the theory can be easily specialized to explicitly accommodate applications in different scenarios. For example, the sub-Gaussian condition would be an option but not required.

In Section 2.1, we study the \textit{conservative bootstrap} where we slightly inflate $t_{\alpha}^*$ by a small fraction, that is, by $1\%$ to $1.01\,t^*_\alpha$. We define an one-sided coverage error in conservative bootstrap as
\begin{align}\label{conservative-bootstrap-error}
\eta^*_{n, \alpha} := \max \Big[ 0,  (1- \alpha) - \mathbb{P} \big\{ T_n \le {1.01}\,t_{\alpha}^* \big\}  \Big]
\end{align}
and derive upper bounds for this quantity.

As $T_n=\max_{j\le p} \sqrt{n}(\overline{X}_{n,j}-\mathbb{E} \overline{X}_{n,j})$, \eqref{conservative-bootstrap-error} means that when a statistician inflates the size of the nominal one-sided confidence band $T_n \le t^*_\alpha$ slightly to $T_n \le 1.01\,t^*_\alpha$ for simultaneous inference about the means $\mathbb{E} \overline{X}_{n,j}, 1 \le j\le p$, the true coverage probability is guaranteed to be at least $1-\alpha -\eta^*_{n, \alpha}$. When the coverage probability is uncertain for the exact bootstrap with $T_n\le t^*_\alpha$, such a small change in the length of confidence interval or the size of rejection region in hypothesis testing would be readily acceptable in real applications when $\eta^*_{n, \alpha}$ is theoretically guaranteed to be small. The thrust of our theoretical result is that in exchange of such a small loss in statistical efficiency, a dramatic reduction in the sample size requirement and guaranteed coverage error would materialize under proper conditions.  

As a matter of fact, we prove that for any $0 < \eta_0 < 1- \alpha$ and $\log (np) \le c_0 n$
\begin{align}\label{intro:conservative_bootstrap}
\eta^*_{n, \alpha} \le C_{c_0, \eta_0} \min \bigg\{ 
& \Big( \frac{(\log (np))^3}{n} \Big)^{1/2} \frac{M^2}{t_{\alpha + \eta_0}^2}, \Big(\frac{(\log (np))^3(\log p)}{n}\Big)^{1/2}\frac{M^2}{t_{\alpha + \eta_0} \overline{\sigma}}, 
\\ \nonumber
& \hspace{12em} \Big(\frac{(\log (np))^3(\log p)^2}{n} \Big)^{1/4} \frac{M}{\overline{\sigma}} \bigg\}
\end{align}
holds for the empirical bootstrap or the multiplier bootstrap with third moment match, where $t_{\alpha + \eta_0}$ is the $(1- \alpha - \eta_0)$-quantile of $T_n$ that should be positive, $M$ depends on 
the moment or tail probability condition on the data and $\overline{\sigma}$ is the soft minimum of the standard deviations $\sigma_j := \big\{\mathrm{Var}(\sqrt{n}\,\overline{X}_{n,j})\big\}^{1/2}$ 
as defined in \eqref{def:soft_minimum} below. 
To express this error bound in full strength and exhibit the rate more clearly, consider the case where certain sub-vectors $(X_{i,j},j\in S)^\top$ satisfy the conditions in \citep{jiang2004asymptotic,xiao2013asymptotic} and $M$ is a constant. If in addition $\log n\lesssim \log|S|\asymp \log p$, then $\log(np) \asymp \log p\lesssim t_{\alpha + \eta_0}$ so that the convergence rate 
\eqref{intro:conservative_bootstrap} is 
$\eta^*_{n, \alpha} \le C \sqrt{(\log p)/n}$.
Even if we are content with constant $t_{\alpha + \eta_0}/\overline{\sigma}$, \eqref{intro:conservative_bootstrap} still yields the error rate $\big( (\log (np))^3 /n \big)^{1/2}$, clearly faster than the 
existing rate for the exact bootstrap. 
More examples are discussed in Section 2.3. 

Note that we inflate the bootstrap quantile $t_{\alpha}^*$ by a fixed 
factor $\epsilon_0 = 0.01$. In fact, our theoretical analyses are carried out for the full spectrum of $\epsilon_0 \in [0, \infty)$, including the case of $\epsilon_0 = \epsilon_n$. Consider the error 
$\eta_{n, \alpha}^*(\epsilon_n) = \max \big[ 0,  (1- \alpha) 
- \mathbb{P} \big\{ T_n \le (1 + \epsilon_n) t_{\alpha}^* \big\}  \big]$. It turns out that our analysis suggests a phase transition in $\eta_{n, \alpha}^*(\epsilon_n)$: If $\epsilon_n$ is of the order $n^{-1/4}$ or smaller, 
our upper bound for the coverage error for the conservative bootstrap does not improve upon the one for the exact bootstrap procedures. For definiteness, we recommend using the \textit{constant} $\epsilon_n = 0.01$ in practice. This is based on the following rationales: (i) The constant is very small 
so that the loss of statistical efficiency is minimal; (ii) The conservative confidence interval is indeed sightly conservative in our simulation experiments for moderate and large $p$; (iii) The small inflation in the size of the confidence band already translates into significant improvements 
in coverage probability -- it provides a healthy cushion of about $1\% - 3\%$ increase in coverage probability in our simulation experiments.


In Section 2.2, we establish as the second main result of this paper improved convergence rates for the exact bootstrap procedures. For the empirical bootstrap and the multiplier bootstrap with third moment match, 
we prove that 
\begin{align}\label{intro:exact_bootstrap}
\big| (1- \alpha) - \mathbb{P} \{ T_n \le t_{\alpha}^* \} \big| \le C \Big(\frac{(\log (np))^3(\log p)^2}{n} \Big)^{1/4} \frac{M}{\overline{\sigma}}, 
\end{align}
in contrast to \eqref{intro:conservative_bootstrap}. 
\cite{koike2019notes} proved that, using a multiplier that is ``close'' to Gaussian but also with third moment match as a `bridge', a triangle inequality can be established to show that the Gaussian multiplier bootstrap enjoys the error bound 
\begin{align}\label{rate-koike}
\Big| \mathbb{P} \big\{ T_n \le t_{\alpha}^* \big\} - (1- \alpha) \Big| \le C \Big(\frac{(\log (np))^5}{n}\Big)^{1/6},
\end{align} 
under suitable moment and tail probability conditions. Later \cite{chernozhukov2019improved} improved the exponent of the above rate from $1/6$ to $1/4$ for general multiplier bootstrap and the empirical bootstrap when the data $\{X_i\}$ are sub-Gaussian. Their results differ from our \eqref{intro:exact_bootstrap} in the following aspects. First, we do not require $\{X_i\}$ to be sub-Gaussian. It is unclear if the results of \cite{chernozhukov2019improved} would hold under more general moment and tail probability conditions as their analysis highly intertwines with the sub-Gaussian condition on $\{X_i\}$. Secondly, even in the sub-Gaussian case, \eqref{intro:exact_bootstrap} is still sharper 
due to the choice of $M$ in \eqref{intro:exact_bootstrap} as discussed in detail in Remark \ref{remark:CCKK} in Section 2.3. 
Moreover, our simulation results in Section 4 suggests that the multiplier bootstrap with third moment match and the empirical bootstrap generally perform better than the Gaussian multiplier bootstrap and the Rademacher multiplier bootstrap; see also \cite{chernozhukov2019improved} and \cite{deng2017beyond} for similar simulation results.

\begin{figure}[h!]
\centering
	\begin{tikzpicture}
	\begin{axis}[
	axis lines=middle,
	ylabel style={anchor = east},
	xlabel=$\varepsilon$,
	xlabel style={anchor = north},
	ymin = 0, ymax = 3.5,
	xmin = 0, xmax = 7.5,
	ytick={1, 3},
	yticklabels={\tiny{$\displaystyle \Big(\frac{(\log (np))^3(\log p)^2}{n} \Big)^{1/2} \frac{M^2}{\overline{\sigma}^2}$}, \tiny{$\displaystyle \Big(\frac{(\log (np))^3(\log p)^2}{n} \Big)^{1/4} \frac{M}{\overline{\sigma}}$}},
	xtick={1, 3, 4.5},
	xticklabels={\tiny{$\displaystyle \Big(\frac{(\log (np))^{3}}{n}\Big)^{1/4} M$}, \tiny{$\displaystyle \frac{\overline{\sigma}}{\sqrt{\log p}}$}, \tiny{$\displaystyle \overline{\sigma}$}},
	x tick label style={anchor = north},
	]
	\addplot [black, line width = 1.5pt, domain=0:1]
	{3};
	\addplot+[ycomb, blue, line width = 1.5pt, opacity = 0.4, mark= *, dashed] plot coordinates
	 {(1,3) (3,1)};
	\addplot+[xcomb, blue, line width = 1.5pt, opacity = 0.4, mark= *, dashed] plot coordinates
	 {(3, 1)};
	\node (source1) at (axis cs:2.5,3) [anchor=west] {\tiny{$\displaystyle \Big(\frac{(\log (np))^3(\log p)}{n}\Big)^{1/2}\frac{M^2}{\varepsilon \overline{\sigma}}$}};
	\node (destination1) at (axis cs:2,1.5){};
	\draw[->, line width = 1.5pt, opacity = 0.5, dashed](source1)--(destination1);
	\node (source2) at (axis cs:4,1.5) [anchor=west] {\tiny{$\displaystyle \Big( \frac{(\log (np))^3}{n} \Big)^{1/2} \frac{M^2}{\varepsilon^2}$}};
	\node (destination2) at (axis cs:5,0.4){};
	\draw[->, line width = 1.5pt, opacity = 0.5, dashed](source2)--(destination2);
	\addplot [black, line width = 1.5pt, domain=1:3]
	{3/x}; 
	\addplot [black, line width = 1.5pt, domain=3:7]
	{9/x^2};
	\end{axis}
	\end{tikzpicture}
	\caption{LP pre-distance upper bound $\overline{\eta}_n(\varepsilon)$.}
	\label{fig:phase_transition}
\end{figure}
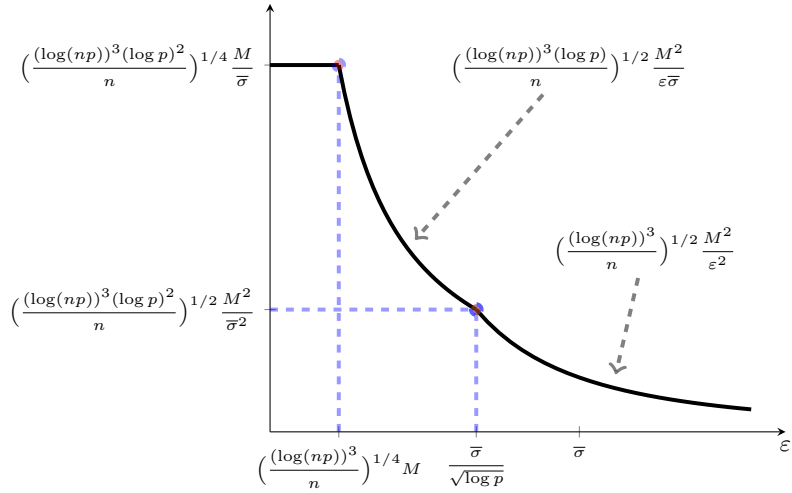

Our main results \eqref{intro:conservative_bootstrap} and \eqref{intro:exact_bootstrap} in this paper are established based on upper bounds on the Levy-Prokhorov (LP) pre-distance between $T_n$ and its bootstrap counterpart $T_n^*$. 
For studying the bootstrap consistency, the LP pre-distance is defined as
\begin{align}\label{intro:LP}
\eta_n^*(\varepsilon, t) &= \max \Big\{0, \mathbb{P} \big\{ T_n \le t -\varepsilon \big\} - \mathbb{P}^* \big\{ T_n^* \le t \big\}, \mathbb{P}^* \big\{ T_n^* \le t -\varepsilon \big\} - \mathbb{P} \big\{ T_n \le t \big\} \Big\}.
\end{align}
This quantity is a key ingredient in the development of our main results as it is closely related to a two-sided version of \eqref{conservative-bootstrap-error} in both appearance and our analysis. Its connection to the consistency of conservative bootstrap in \eqref{intro:conservative_bootstrap} and that of exact bootstrap in \eqref{intro:exact_bootstrap} are discussed in detail in Section 2.4. We prove that with high probability $\eta_n^*(\varepsilon, t) \le C\overline{\eta}_n(\varepsilon)$ where
\begin{align*}
\overline{\eta}_n(\varepsilon) :=  \min \bigg\{ & \Big( \frac{(\log (np))^3}{n} \Big)^{1/2} \frac{M^2}{\varepsilon^2}, \Big(\frac{(\log (np))^3(\log p)}{n}\Big)^{1/2}\frac{M^2}{\varepsilon \overline{\sigma}}, \Big(\frac{(\log (np))^3(\log p)^2}{n} \Big)^{1/4} \frac{M}{\overline{\sigma}} \bigg\},
\end{align*}
which is represented by the three-piece curve in Figure \ref{fig:phase_transition}. It improves upon the larger bound $\overline{\eta}^{(DZ)}_n(\varepsilon) = \big({(\log (np))^3}/{n} \big)^{1/2} {M^2}/{\varepsilon^2}$ derived in \cite{deng2017beyond} when $\varepsilon \le \overline{\sigma}/\sqrt{\log p}$. 

We bound the LP pre-distance $\eta_n^*(\varepsilon, t)$ through a comparison theorem which provides upper bounds for the absolute difference between the expectations of functions of $T_n$ and its bootstrap analogue $T_n^*$, 
\begin{align}\label{intro:bound_LP}
\eta_n^*(\varepsilon, t) \le \big| \mathbb{E} h_t(T_n) - \mathbb{E}^* h_t(T_n^*) \big|,
\end{align}
where $h_t(x)$ is a smooth function with bounded derivatives up to a certain order that coincides with the indicator function $I\{x \le t\}$ outside a small interval $[t - \varepsilon, t]$. The problem of bounding the LP pre-distance then boils down to the derivation of a suitable comparison theorem. To this end, \cite{deng2017beyond} developed the coherent Lindeberg interpolation method for deriving the comparison bound for independent matrices $\{(X_i, Y_i) \in \mathbb{R}^{p\times 2}, 1\le i\le n\}$, that is, $ \mathbb{E} f(X_1, \ldots, X_n) -  \mathbb{E} f(Y_1, \ldots, Y_n)$ where $f$ is a smooth function. 
 This coherent Lindeberg interpolation has been shown in \cite{deng2017beyond} to be a powerful tool in developing comparison bounds. However, their analysis can still be improved by noticing that the derivatives of $h_t$ are all zero except in the small region $[t - \varepsilon, t]$. More precisely, the $m$-th derivative of $h_t$ can be written as
\begin{align}\label{intro:comparison_tight}
h_t^{(m)}(x) = h_t^{(m)}(x) I\{ t - \varepsilon \le x \le t\}, \forall m \ge 1
\end{align}
which will potentially tighten the upper bound after we apply interpolation arguments and Taylor series expansion on \eqref{intro:bound_LP}. This observation is inspired by \cite{bentkus2003dependence} who took advantage of a similar observation and carried out an exquisite analysis leading to a nearly optimal Berry-Esseen bound for Gaussian approximation on convex sets. The idea is also utilized in \cite{chernozhukov2017central}. 


In Section 3, we throughly study the incorporation of the above \cite{bentkus2003dependence} idea to the coherent Lindeberg interpolation of \cite{deng2017beyond} to establish various comparison theorems.
An important outcome of this study is a new permutation invariance lemma, Lemma \ref{lm:PI} in Section 3.1, which produces more general comparison bounds, Theorems \ref{thm:weak_comparison} and \ref{thm:strong_comparison}, in Section 3.2 to allow broader applications of the coherent Lindeberg interpolation, including utilizing the important argument \eqref{intro:comparison_tight}. The improved coherent interpolation method contributes to both of our main results, \eqref{intro:conservative_bootstrap} and \eqref{intro:exact_bootstrap}, and will be a very useful tool for deriving more comparison bounds in other problems. 

The rest of the paper is organized as follows. We give the main results of this paper, the consistency of conservative and exact bootstrap, in Section 2. Examples are also given in this section to derive more explicit rates. In Section 3, we give a detailed discussion on the improved coherent Lindeberg interpolation method and develop several crucial comparison theorems towards our main results. Some numerical experiments are given in Section 4 to support our theoretical results. The technical proofs are relegated to the Appendix.

\textit{Notation.} We use the following notation.
We assume $n \to \infty$ and $p=p_n$ to allow $p \to \infty$ as $n \to \infty$. We assume $\log p \ge 1$ throughout the paper; it can be replaced with $(\log p) \vee 1$ to allow $\log p < 1$. For any vector $x = (x_1, \ldots, x_p)^\top \in \mathbb{R}^p$, $x^m = (x_1^m, \ldots, x_p^m)^\top$ is still a $p$ dimensional vector, and $x^{\otimes m} = (x_{i_1}\cdots x_{i_m})_{p \times \cdots \times p}$ is a tensor of order $m$. We write partial derivative operators as tensors 
$(\partial/\partial x)^{\otimes m}=\big((\partial/\partial x_{i_1})\cdots (\partial/\partial x_{i_m})\big)_{p\times \cdots \times p}$ 
for $x=(x_1,\ldots,x_p)^\top$, 
so that ${f}^{(m)} = (\partial/\partial x)^{\otimes m}f(x)$ is a tensor for functions $f(x)$ of input $x\in\mathbb{R}^p$,
and for two $m$-th order tensors $f$ and $g$ in $\mathbb{R}^{p\times \cdots \times p}$,
the vectorized inner product is denoted by
$\big\langle f,g \big\rangle = \sum_{j_1 =1}^p \cdots \sum_{j_m=1}^p  f_{j_1,\ldots, j_m} g_{j_1,\ldots, j_m}$
and $|f|\le |g|$ means $|f_{j_1,\ldots, j_m}|\le |g_{j_1,\ldots, j_m}|$ for all indices ${j_1,\ldots, j_m}$. 
Let $\|\cdot\|_{p}$ be the $\ell_p$ norm and $\|\cdot\|_{\infty}$ the maximum norm in finite dimensional space. For tensors, we may vectorize them and apply these norms. The vectorized maximum norm is denoted by $\| \cdot \|_{\max}$. We denote by $C$ a numerical constant and $C_{\rm index}$ a constant depending on the ``index'' only. For example, $C_{a,b,c}$ is a constant depending on $(a,b,c)$ only.
Finally, we define the maximum average moment of $\{X_i\}$ as 
\begin{align}\label{average_moments}
\textstyle M_m^m := \max_{1\le j \le p} n^{-1} \sum_{i=1}^n \mathbb{E} |X_{i,j} {- \mathbb{E} X_{i,j}}|^m, 
\end{align}
Let $\sigma_j^2 = n^{-1} \sum_{i=1}^n \mathrm{Var}(X_i)$, and denote the $j$-th smallest standard deviation by $\sigma_{(j)}$ and the minimum by $\underline{\sigma} = \sigma_{(1)}$. We define the soft minimum as
\begin{align}\label{def:soft_minimum}
\overline \sigma = \min_{1\le j \le p} \Big\{ \big( 2 + \sqrt{2 \log p} \big) \big/ \big[ 1/\sigma_{(1)} +  (1 + \sqrt{2\log j})/\sigma_{(j)} \big] \Big\}
\end{align}
which satisfies $\overline{\sigma} \ge \underline{\sigma}$.

\section{Consistency of conservative and exact bootstrap}
We study several different bootstrap schemes, mainly the empirical bootstrap and the multiplier bootstrap with different multipliers. Recall that $\overline{X}_n = \sum_{i=1}^n X_i/n$.

\textit{(i).\,Efron's empirical bootstrap} \citep{efron1979}. We generate i.i.d.\,vectors $X_1^*,\ldots,X_n^*$ from the empirical distribution of the centered sample $X_1- \overline{X}_n,\ldots,X_n - \overline{X}_n$: Under the conditional probability $\mathbb{P}^* = \mathbb{P}\{ \ \cdot \ | X_i \, \forall i\}$,
\begin{align}\label{empirical_bootstrap}
\textstyle \mathbb{P}^* \big\{X_i^* = X_k - \overline{X}_n  \big\} = n^{-1} \sum_{\ell=1}^n I\{X_\ell = X_k\} , \ 1 \le k, i \le n. 
\end{align}

\smallskip
\textit{(ii).\,Multiplier/Wild bootstrap}. In multiplier/wild bootstrap \citep{wu1986jackknife}, we generate
\begin{align}\label{multiplier_bootstrap}
X_i^* = W_i\big(X_i - \overline{X}_n \big), \ i=1, \ldots, n,
\end{align}
where $W_1,\ldots, W_n$ are i.i.d.~random variables with
$\mathbb{E} \,W_i = 0, \mathbb{E} \,W_i^2=1$,
and the sequence $\{W_i\}$ is independent of the original data $\{X_i\}$. Based on the choice of the multipliers, we list below some common bootstrap schemes:
\begin{itemize}
		\item \textit{Gaussian multiplier bootstrap}, where the multipliers are standard Gaussian variables, that is, $W_i \sim \mathcal{N}(0,1)$ for all $i$\,;
		\item \textit{Rademacher's multiplier bootstrap}, where the multipliers are Rademacher variables, that is, $ \mathbb{P} \{W_i = \pm 1\} = 1/2 $ for all $i$\,;
		\item \textit{Mammen's multiplier bootstrap}, where the multiplier satisfies
		$\mathbb{P} \big\{W_i= (1\pm\sqrt{5})/2 \big\}= (\sqrt{5}\mp1)/(2\sqrt{5})$ for all $i$.
		This multiplier, proposed by \cite{mammenwild1993}, satisfies the third moment match condition \citep{liu1988bootstrap}, that is, $\mathbb{E} \, W_i^3 =1$.
	\end{itemize}

According to their moment match behavior, we hereafter refer to any multiplier bootstrap with sub-Gaussian multipliers $W_i$ satisfying 
\begin{align}\label{third_moment_match}
\sum_{i=1}^n \mathbb{E} W_i^3(X_i - \mathbb{E} X_i )^{\otimes 3}/n = \sum_{i=1}^n \mathbb{E} (X_i - \mathbb{E} X_i )^{\otimes 3}/n
\end{align}
the \textit{multiplier bootstrap with third moment match}, including the Mammen's multiplier bootstrap and any other multiplier bootstrap with sub-Gaussian multipliers satisfying $\mathbb{E} W_i^3 =1, \forall i$. The Gaussian multiplier bootstrap and the Rademacher's multiplier bootstrap are also multiplier bootstrap schemes with third moment match when $\{X_i\}$ have the average third moment zero, that is, $\sum_{i=1}^n \mathbb{E} ( X_i - \mathbb{E} X_i)^{\otimes 3}/n = 0$. 

Corresponding to the maximum of the sums, $T_n$, we define its bootstrap analogue as 
$T_{n}^* = \max_j \sum_{i=1}^n X^*_{i,j} /\sqrt{n}$,
where $\{X_i^*\}$ is a bootstrap sample drawn from the bootstrap schemes mentioned above. Let $t_{\alpha} := \inf \big[ t: \mathbb{P} \{T_n > t\} \le \alpha\big]$ be the $(1- \alpha)$-quantile of $T_n$ and similarly $t_{\alpha}^*$ of $T_n^*$. 

\subsection{Consistency of conservative bootstrap} 

In conservative bootstrap, we slightly inflate $t_{\alpha}^*$ and use $(1+ \epsilon_0) t_{\alpha}^* = 1.01 t_{\alpha}^*$. Here $\epsilon_0 = 0.01$ is called the inflation factor. We are interested in the one-sided coverage error of conservative bootstrap
$\eta^*_{n, \alpha} = \max \big[ 0,  (1- \alpha) - \mathbb{P} \{ T_n \le (1 + 0.01) t_{\alpha}^*\}  \big]$ as in \eqref{conservative-bootstrap-error}.

Let constant $C_{\mathrm{boot}}$ be a numerical constant $C_0$ for empirical bootstrap and a constant $C_{\tau_0}$ for multiplier bootstrap where i.i.d.\,multipliers $\{W_i\}$ satisfy the sub-Gaussian condition $\mathbb{E} \exp(tW_1) \le \exp(\tau_0^2 t^2/2) \ \forall t$. We have the following theorem.

\begin{theorem}\label{thm:main}
Let $\{ X_i\in\mathbb{R}^p, 1\le i \le n\}$ be independent random vectors and $\{ X_i^*\}$ generated by the empirical bootstrap or the multiplier bootstrap with third moment match. Assume $\log(np) \le c_0\,n$ for a fixed constant $c_0$. Recall $t_{\alpha}$ is the $(1- \alpha)$ quantile of $T_n$. For $\varepsilon > 0$ and $M \ge M_4 := \max_j \sum_{i=1}^n \mathbb{E} |X_{i,j} - \mathbb{E} X_{i,j}|^4/n$, define
\begin{align}\label{definition-prob-thm:main}
q_0(M, \varepsilon) =  \mathbb{P} \Big\{\| \boldsymbol{X}  - \mathbb{E}  \boldsymbol{X} \|_{\max}> \max\Big[M\big(n/\log (np) \big)^{1/4}, \sqrt{n} \varepsilon \big/\log (np) \Big]\Big\}.
\end{align}
Then, for any $0 < \eta_0 < 1 - \alpha$ such that $t_{\alpha + \eta_0} > 0$
\begin{align}\label{bound-thm:main}
\eta^*_{n, \alpha} 
&= \max \Big[ 0,  (1- \alpha) - \mathbb{P} \big\{ T_n \le 1.01 t_{\alpha}^* \big\}  \Big] 
\\ \nonumber
&\le C_{\mathrm{boot}, c_0, \eta_0} \Bigg( \min \bigg\{ \Big( \frac{(\log (np))^3}{n} \Big)^{1/2} \frac{M^2}{t_{\alpha + \eta_0}^2}, \Big(\frac{(\log (np))^3(\log p)}{n}\Big)^{1/2}\frac{M^2}{t_{\alpha + \eta_0} \overline{\sigma}},
\\ \nonumber
& \hspace{8em} \Big(\frac{(\log (np))^3(\log p)^2}{n} \Big)^{1/4} \frac{M}{\overline{\sigma}} \bigg\} + \frac{1}{np} + q_0\big(M, t_{\alpha + \eta_0}/101 \big) \Bigg),
\end{align}
where $C_{\mathrm{boot}, c_0, \eta_0}$ depends on the bootstrap, $c_0$ and $\eta_0$. The soft minimum standard deviation $\overline{\sigma}$ is defined as in \eqref{def:soft_minimum}.
\end{theorem}

\begin{remark}\label{remark:main}
We can explicitly write the minimization term on the right-hand side of \eqref{bound-thm:main} as
\begin{align*}
\begin{cases}
\displaystyle \Big( \frac{(\log (np))^3}{n} \Big)^{1/2} \frac{M^2}{t_{\alpha + \eta_0}^2}, & \displaystyle t_{\alpha + \eta_0} \ge \frac{\overline{\sigma}}{\sqrt{\log p}},
\\[10pt]
\displaystyle \Big(\frac{(\log (np))^3(\log p)}{n}\Big)^{1/2}\frac{M^2}{t_{\alpha + \eta_0} \overline{\sigma}}, & \displaystyle \frac{(\log (np))^{3/4}}{n^{1/4}}M \le t_{\alpha + \eta_0} \le  \frac{\overline{\sigma}}{\sqrt{\log p}},
\\[10pt]
\displaystyle \Big(\frac{(\log (np))^3(\log p)^2}{n} \Big)^{1/4} \frac{M}{\overline{\sigma}}, & \displaystyle t_{\alpha + \eta_0}
 \le \frac{(\log (np))^{3/4}}{n^{1/4}}M .
\end{cases}
\end{align*}
which is $\overline{\eta}_n(t_{\alpha + \eta_0})$ in Figure \ref{fig:phase_transition}.
\end{remark}

Note that the $(1- \alpha - \eta_0)$-quantile $t_{\alpha + \eta_0}$ is much likely to grow with $p$, slow or fast. Even if we take a pessimistic lower bound of $t_{\alpha + \eta_0}$ as roughly of order $\overline{\sigma}$ and assume $M$ is bounded, Theorem \ref{thm:main} asserts that 
the conservative bootstrap procedure has error bound $\eta_{n, \alpha}^* \le \big( (\log (np))^3 /n \big)^{1/2}$.
When $t_{\alpha + \eta_0} \asymp \sqrt{\log p} \, \overline{\sigma}$ and $\log p \asymp \log(np)$, it becomes
$\eta_{n, \alpha}^* \le ((\log p)/n)^{1/2}$.
In Section 2.3, we impose different conditions on $\{X_i\}$ and further specify the error bounds based on the comprehensive Theorem \ref{thm:main}. 

To the best of our knowledge, the above results for conservative bootstrap are new. It is tempting to wonder if we can carry out a conservative inference by adding a small amount to the bootstrap quantile as $t_{\alpha}^* + \varepsilon$, instead of inflating $t_{\alpha}^*$ as we proposed here. \cite{chernozhukov2019improved} studied a similar procedure where the shift $\varepsilon$, referred to as the infinitesimal factor, is a constant and $\{X_i\}$ are sub-Gaussian. However, the approach requires the knowledge of the scale of $t_{\alpha}^*$, which is not usually feasible in practice. 

As mentioned in the Introduction Section, the inflation factor $\epsilon_0$ can be arbitrary in theory. A even more comprehensive version of Theorem \ref{thm:main} is given below.

\begin{theorem}\label{thm:main-full}
For any $\epsilon_n \ge 0$, define
$\eta^*_{n, \alpha}(\epsilon_n) = \max \Big[ 0,  (1- \alpha) - \mathbb{P} \big\{ T_n \le (1 + \epsilon_n) t_{\alpha}^* \big\}  \Big]$.
The bound \eqref{bound-thm:main} in Theorem \ref{thm:main} holds for $\eta^*_{n, \alpha} = \eta^*_{n, \alpha}(\epsilon_n)$ and $t_{\alpha + \eta_0}/101$ replaced with $\epsilon_n t_{\alpha + \eta_0}/(1 + \epsilon_n)$.
\end{theorem}


In view of Remark \ref{remark:main}, we know the bound given in Theorem \ref{thm:main-full} for $\eta_{n, \alpha}^*(\epsilon_0)$ remains a three-piece function. It is worth noting that for really small $\epsilon_n$, that is, when $\epsilon_n/(1 + \epsilon_n) \le \{(\log (np))^3/n\}^{1/4}(M/t_{\alpha + \eta_0})$, the bound becomes the worst possible one that is essentially the two-sided coverage error bound of exact bootstrap, given in the forthcoming subsection.


\subsection{Consistency of exact bootstrap}
In exact bootstrap, we directly use the bootstrap quantile $t_{\alpha}^*$ for inference and consider the two-sided coverage error. We have the following theorem.

\begin{theorem}\label{thm:main-exact}
Let $\{ X_i\in\mathbb{R}^p, 1\le i \le n\}$ be independent random vectors and $\{ X_i^*\}$ generated by the empirical bootstrap or the multiplier bootstrap with third moment match. For any $M \ge M_4$ as in \eqref{average_moments}, there exists a constant $C_{\mathrm{boot}}$ depending on the bootstrap as in Theorem \ref{thm:main} such that
\begin{align}\label{bound-thm:main-exact}
& \Big| (1- \alpha) - \mathbb{P} \big\{ T_n \le t_{\alpha}^* \big\}  \Big| 
\\ \nonumber
&\le 4 \,  \mathbb{P} \Big\{\| \boldsymbol{X}  - \mathbb{E}  \boldsymbol{X} \|_{\max}> M\big( n / \log (np) \big)^{1/4}\Big\} + C_{\mathrm{boot}} \Big(\frac{(\log (np))^3(\log p)^2}{n} \Big)^{1/4} \frac{M}{\overline{\sigma}} .
\end{align}
\end{theorem}


\begin{remark}
\cite{deng2017beyond} gives a universal upper bound 
\begin{align*}
\Big( \frac{(\log (np))^5}{n} \Big)^{1/6} \frac{ \big( n^{-1} \sum_{i=1}^n \mathbb{E} \max_{1\le j \le p} |X_{i,j} - \mathbb{E} X_{i,j}|^4 \big)^{1/4}}{\overline{\sigma}}
\end{align*}
for the empirical bootstrap and another smaller universal bound 
\begin{align*}
\Big( \frac{(\log (np))^5}{n} \Big)^{1/6} \frac{ \big( \mathbb{E} \max_{1\le j \le p} n^{-1} \sum_{i=1}^n |X_{i,j} - \mathbb{E} X_{i,j}|^4 \big)^{1/4}}{\overline{\sigma}}
\end{align*} 
for any multiplier bootstrap with third moment match and sub-Gaussian multipliers. They become effective when $\| \boldsymbol{X}  - \mathbb{E}  \boldsymbol{X} \|_{\max}$ has a really heavy tail, that is, when the resulting $M$ in \eqref{bound-thm:main-exact} is too large. These bounds, however, seems hard to improve using our new proof technique. We may take the minimum of \eqref{bound-thm:main-exact} and the above universal bounds as the final upper bounds for the consistency of exact bootstrap.
\end{remark}

\subsection{Examples}\label{subsection:example}
We give some examples in this subsection to derive explicit consistency rates for $\{X_i\}$ under various conditions. To apply Theorems \ref{thm:main} and \ref{thm:main-exact} in full strength, we then consider the problem of approximating the distribution of the maximum deviation of sample covariance estimates with a specified order of $t_{\alpha + \eta_0}$.

\begin{example}
In this example, $\overline{\sigma}$ and $M_4$ are allowed to depend on $n$ and to diverge to $0$ or $\infty$, but they can also be treated as constants for simplicity.
We consider four examples specified by certain measure $B_n$ of the tail of $\{X_{i,j}\}$:

\begin{enumerate}[label=(E.{\arabic*}),leftmargin=0.8in, 
itemsep=6pt
]
	\item $\mathbb{P} \big\{ |X_{i,j} - \mathbb{E} X_{i,j}| \le B_n \big\} =1$ \label{ex-bounded};

	\item $\mathbb{E} \exp \big[ |X_{i,j} - \mathbb{E} X_{i,j}|^2/B_n^2 \big] \le 2 $ for all $i,j$; \label{ex-gaussian}

	\item $\mathbb{E} \exp \big[ |X_{i,j} - \mathbb{E} X_{i,j}|/B_n \big] \le 2 $ for all $i,j$; \label{ex-exp}

	\item $n^{-1} \sum_{i=1}^n \mathbb{E} \, \max_{1\le j \le p} |X_{i,j} -\mathbb{E} X_{i,j}|^q \le B_n^q$. \label{ex-moment}
\end{enumerate}

Due to the three-piece structure of the bound in \eqref{bound-thm:main}, the complete upper bounds for $\{X_i\}$ satisfying either one condition among \ref{ex-bounded}-\ref{ex-moment} can be quite complicated. To capture the essence of the advantage of conservative bootstrap procedure, and since $t_{\alpha + \eta_0} \gtrsim \overline{\sigma}/\sqrt{\log p}$ is most likely to happen, we only consider the first component in the minimization term of \eqref{bound-thm:main}, $\{ (\log (np))^3/n \}^{1/2} (M^2/t_{\alpha + \eta_0}^2)$, in the following corollary. 

\begin{corollary}\label{coro:example-1}
Let $\{ X_i\in\mathbb{R}^p, 1\le i \le n\}$ be independent random vectors and $\{ X_i^* \}$ generated by the empirical bootstrap or the multiplier bootstrap with third moment match. Let $a \lesssim_{boot} b$ represents $a \le C_{\mathrm{boot}} b$ where $C_{\mathrm{boot}}$ depends on the bootstrap as in Theorem \ref{thm:main}. The one-sided coverage error of conservative bootstrap is bounded by 
\begin{align*}
\eta_{n, \alpha}^* &= \max \big[ 0,  (1- \alpha) - \mathbb{P} \{ T_n \le 1.01 t_{\alpha}^* \}  \big]
\\ \nonumber
&\lesssim_{\mathrm{boot}, c_0, \eta_0} \Big( \frac{ (\log (np))^3}{n}  \Big)^{1/2} \frac{M_4^2}{t_{\alpha + \eta_0}^2} + 
\begin{cases} 
\big\{ (\log (np))^2 / n \big\} (B_n^2/ t_{\alpha + \eta_0}^2), & \hbox{ under \ref{ex-bounded}}  
\\[6pt]
\big\{ (\log (np))^3 / n \big\} (B_n^2/ t_{\alpha + \eta_0}^2), & \hbox{ under \ref{ex-gaussian}} 
\\[6pt]
\big\{ (\log (np))^4 / n \big\} (B_n^2/ t_{\alpha + \eta_0}^2), & \hbox{ under \ref{ex-exp}} 
\\[6pt]
C_q \big\{(\log(np))^q / n^{q/2-1} \big\} (B_n^q/ t_{\alpha + \eta_0}^q), & \hbox{ under \ref{ex-moment}}
\end{cases}
\end{align*}
where we assume $\log(np) \le c_0n$.
The two-sided coverage error of exact bootstrap is bounded by
\begin{align*}
& \big| \mathbb{P} \{ T_n \le t_{\alpha}^* \} - (1- \alpha) \big| 
\\ \nonumber
& \lesssim_{\mathrm{boot}} \Big(\frac{(\log (np))^3(\log p)^2}{n} \Big)^{1/4} \frac{M_4}{\overline{\sigma}} + \begin{cases} 
\big\{(\log (np))^2(\log p)/ n \big\}^{1/2} (B_n /\overline{\sigma}), & \hbox{ under \ref{ex-bounded}}  
\\[6pt]
\big\{(\log (np))^3(\log p)/ n \big\}^{1/2} (B_n /\overline{\sigma}) , & \hbox{ under \ref{ex-gaussian}} 
\\[6pt]
\big\{(\log (np))^4(\log p)/ n \big\}^{1/2} (B_n /\overline{\sigma}) , & \hbox{ under \ref{ex-exp}} 
\\[6pt]
\Big\{ \big[ (\log(np)) (\log p)^{1/2} / n^{1/2-1/q} \big] (B_n /\overline{\sigma}) \Big\}^{q/(q+1)}. & \hbox{ under \ref{ex-moment}}
\end{cases}
\end{align*}
\end{corollary}

\end{example}

\begin{remark}\label{remark:CCKK}
\cite{chernozhukov2019improved} studies the consistency of exact bootstrap under sub-Gaussian condition \ref{ex-gaussian}. Their equivalent statement is that, for empirical bootstrap and multiplier bootstrap, $\big| \mathbb{P} \{ T_n \le t_{\alpha}^* \} - (1- \alpha) \big| \le C_{\mathrm{boot}} \, \gamma$ holds with
\begin{align*}
\gamma = \gamma^{(CCKK)} := \Big(\frac{(\log (np))^5}{n} \Big)^{1/4} \frac{M_4}{\underline{\sigma}} + \bigg[ \Big(\frac{(\log (np))^5}{n} \Big)^{1/2} \frac{B_n}{\underline{\sigma}} \bigg]^{1/2},
\end{align*}
where $\underline{\sigma}^2 = \min_j n^{-1}\sum_{i=1}^n \mathrm{Var}(X_{i,j}) \le \overline{\sigma}^2$. This should be compared with the smaller
\begin{align*}
\gamma = \gamma^* := \Big(\frac{(\log (np))^3(\log p)^2}{n} \Big)^{1/4} \frac{M_4}{\overline{\sigma}} + \Big(\frac{(\log (np))^3(\log p)}{n} \Big)^{1/2} \frac{B_n}{\overline{\sigma}} 
\end{align*}
from the above Corollary \ref{coro:example-1}. Despite the improvement from the minimum $\underline{\sigma}$ to the soft minimum $\overline{\sigma}$ in $\gamma^*$, the consistency rate $\gamma^{(CCKK)}$ from \citep{chernozhukov2019improved} may still be slower under some circumstances. For example, when $B_n \asymp n^{1/4} M_4$ and $M_4/\underline{\sigma} = M_4/\overline{\sigma} =$ const., we have
\begin{align*}
\gamma^{(CCKK)} \asymp \Big(\frac{(\log (np))^{10}}{n} \Big)^{1/8} \gtrsim \Big(\frac{(\log (np))^{6} (\log p)^2}{n} \Big)^{1/4} \asymp \gamma^{*}.
\end{align*}
In this case, not only the size requirement is reduced from $ n \gg (\log p)^{10}$ to $n \gg (\log p)^8$ but also the overall consistency rate is much faster with exponent $1/4$. We note that the rate $\gamma^{(CCKK)}$ holds for the empirical bootstrap and \textit{any} multiplier bootstrap, but $\gamma^{*}$ is only valid for the empirical bootstrap and the multiplier bootstrap with third moment match.
\end{remark}

In what follows, we study an specific model and specify the order of $t_{\alpha + \eta_0}$ to show the full strength of Theorem \ref{thm:main}.
Suppose $\boldsymbol{Y} = (Y_{i,j})_{1\le i \le n, 1\le j \le m}$ is a matrix with $n$ i.i.d.\,rows. The sample covariance is $\widehat{\Sigma}_{m \times m}$ with the $(j,k)$-th entry as $\widehat{\sigma}_{j,k} = n^{-1} \sum_{i=1}^n (Y_{i,j} - \overline{Y}_{n,j}) (Y_{i,k} - \overline{Y}_{n,k})$, where $\overline{Y}_n = \sum_{i=1}^n Y_i /n$. We are interested in the distribution of $\max_{1 \le j < k \le m} |\widehat{\sigma}_{j,k} - \sigma_{j,k}|$, where $\sigma_{j,k} = \mathrm{Cov}(Y_{1,j}, Y_{1,k})$ is the true covariance. To simplify our discussion, we ignore negligible terms and consider $\widetilde{M}_n = \max_{1 \le j < k \le m} \Big|n^{-1}\sum_{i=1}^n ( Y_{i,j}Y_{i,k} - \mathbb{E} Y_{i,j} Y_{i,k}) \Big|$.
We vectorize $\{Y_{i,j}Y_{i,k} \hbox{ and } -Y_{i,j}Y_{i,k}, \  1 \le j < k \le m\}$ to $X_i = (X_{i,1}, \ldots, X_{i,p})^\top$ where $p = m(m-1)$, and therefore $\widehat{M}_n$ can be written as the maximum of sums. Recall that $X_1, \ldots, X_n$ are i.i.d.\, so that $\mathbb{E} X_{i,\ell} = \sum_{i=1}^n \mathbb{E} X_{i, \ell}/n$, we have
\begin{align*}
\textstyle \mathbb{P} \big\{ \widetilde{M}_n \le t \big\} = \mathbb{P} \big\{ T_n/\sqrt{n} := \max_{1 \le \ell \le p} n^{-1} \sum_{i=1}^n (X_{i, \ell} - \mathbb{E} X_{i, \ell})  \le t \big\}. 
\end{align*}
We also assume constant variances of all entries, that is, $\sigma_0^2 = \frac{1}{n}\sum_{i=1}^n \mathrm{Var}(Z_{i, \ell})$ for all $1 \le \ell \le p$. It follows from \cite{xiao2013asymptotic} that under the conditions assumed therein the $(1- \alpha - \eta_0)$-quantile of $\sqrt{n} \widetilde{M}_n$ has asymptotic order $\sigma_0 \log m$, or $\sigma_0 \log p$. Suppose $\log (np) \asymp \log p$.
The one-sided coverage error of conservative bootstrap for this $T_n = \sqrt{n}\widetilde{M}_n$ is bounded by
\begin{align*}
\eta_{n, \alpha}^* & = \max \Big[ 0,  (1- \alpha) - \mathbb{P} \big\{ T_n = \sqrt{n} \widetilde{M}_n \le 1.01 t_{\alpha}^* \big\}  \Big]
\\ \nonumber
& 
\lesssim_{\mathrm{boot}} \Big( \frac{\log p}{n} \Big)^{1/2} \frac{M_4^2}{\sigma_0^2} + 
\begin{cases} 
\big\{ (\log  p) / n \big\} (B_n^2/ \sigma_0^2), & \hbox{ under \ref{ex-bounded}}  
\\[6pt]
\big\{ (\log  p)^2 / n \big\} (B_n^2/ \sigma_0^2), & \hbox{ under \ref{ex-gaussian}} 
\\[6pt]
\big\{ (\log  p)^3 / n \big\} (B_n^2/ \sigma_0^2), & \hbox{ under \ref{ex-exp}} 
\\[6pt]
C_q \big\{(\log p)^{q/2} / n^{q/2-1} \big\} (B_n^q/ \sigma_0^q). & \hbox{ under \ref{ex-moment}}
\end{cases}
\end{align*}
When $M_4/\sigma_0$ and $B_n/\sigma_0$ are constants, the sample size requirement for $ \eta_{n, \alpha}^* = o(1)$ is
$(\log p)^{\kappa} \ll n$
with $\kappa =1$ under \ref{ex-bounded}, $\kappa=2$ under \ref{ex-gaussian}, $\kappa = 3$ under \ref{ex-exp} and $\kappa = q/(q-2)$ for $q > 2$ under \ref{ex-moment}. 

\subsection{Connection between consistency of bootstrap and LP pre-distance}
We show the one-sided and two-sided coverage errors in \eqref{conservative-bootstrap-error} and \eqref{bound-thm:main-exact} are closely related to the LP pre-distance in bootstrap
\begin{align}\label{def:bootstrap_LP}
\eta_n^*(\varepsilon, t) &= \max \Big\{0, \mathbb{P} \big\{ T_n \le t -\varepsilon \big\} - \mathbb{P}^* \big\{ T_n^* \le t \big\}, \mathbb{P}^* \big\{ T_n^* \le t -\varepsilon \big\} - \mathbb{P} \big\{ T_n \le t \big\} \Big\}, \ \varepsilon \ge 0.
\end{align}
The following Lemma connects this bootstrap LP pre-distance to the consistency of conservative bootstrap and exact bootstrap.

\begin{lemma}\label{lm:lp-to-consistency}
For $\eta, \epsilon, \varepsilon >0$, the following inequalities hold:
\begin{enumerate}[label = \emph{(\roman*).}, leftmargin = 0.7in, itemsep=1.6ex]
	\item $1-\alpha - \mathbb{P} \big\{T_n\le (1+  \epsilon) t^*_\alpha \big\} \le \mathbb{P} \big\{ \eta_n^* \big(  \epsilon\, t_{\alpha + \eta} / (1+ \epsilon), t_{\alpha + \eta}\big) \ge \eta \big\} + \eta$;
	\item $\big| \mathbb{P} \{T_n\le t^*_\alpha \} - (1-\alpha) \big|  \le \sup_t  \mathbb{P} \big\{ \eta_n^*(\varepsilon, t) > \eta \big\} +\eta + \omega_n(\varepsilon; T_n)$;
	\item $1-\alpha - \mathbb{P} \big\{T_n\le t^*_\alpha +  \varepsilon \big\} 
		\le  \mathbb{P} \big\{ \eta_n^*( \varepsilon, t_{\alpha + \eta}) \ge \eta \big\} + \eta$.
\end{enumerate}
\end{lemma}


As discussed below Theorem \ref{thm:main}, the scale of $t_{\alpha}^*$ is likely unknown and the conservative bootstrap procedure that adds a small amount $\varepsilon$ to $t_{\alpha}^*$ is more difficult to implement in practice. However, (iii) is given here for readers to use in some special cases where this approach is preferred. As this procedure is beyond the scope of this paper and actually its theoretical results are directly visible by this Lemma and Proposition \ref{prop:LP} in Appendix, we omit the details.

With the help of Lemma \ref{lm:lp-to-consistency}, our core task is to bound the bootstrap LP pre-distance $\eta_n^*(\varepsilon, t)$, so that Theorems \ref{thm:main} and \ref{thm:main-full} immediately follow from (i) and Theorem \ref{thm:main-exact} follows from (ii). Notice that
$\eta_n^*(\varepsilon, t) \le \big| \mathbb{E} h_t(T_n) - \mathbb{E}^* h_t(T_n^*) \big|$ as in \eqref{intro:bound_LP}
and $T_n$ can be approximated by the smooth softmax function to be defined later. 
This leads us to consider the general comparison bound of the absolute difference between expectations of smooth functions of $\{X_i\}$ and another data set, $\{Y_i\}$ say, which is treated in the forthcoming Section 3.

\section{Comparison Theory}
Let $(X_i, Y_i)\in \mathbb{R}^{p\times 2}$, $1\le i\le n$ be independent random matrices under $\mathbb{E}$. With a slight abuse of notation, we let $\sigma$ be a permutation operator of $\{1, \ldots, n\}$ in this section, that is, $\{\sigma_1, \ldots, \sigma_n\} = \{1, \ldots, n\}$. A function $f(x_1, \ldots, x_n)$ with $x_i \in \mathbb{R}^p$ is called permutation invariant if $f(x_1, \ldots, x_n) = f(x_{\sigma_1}, \ldots, x_{\sigma_n})$ for any permutation $\sigma$. 

The main task in this section is to derive comparison bounds for
\begin{align*}
\big| \mathbb{E} f(X_1, \ldots, X_n) - \mathbb{E} f(Y_1, \ldots, Y_n) \big|
\end{align*}
and its more general forms, where $f$ is a smooth permutation invariant function. Throughout the entire section, assume $(X_i, Y_i)$ mean zero for all $i$. In Section 3.1, we describe the improved coherent Lindeberg interpolation method whose original version is introduced by \cite{deng2017beyond}. Using the method, two general comparison theorems are derived in Section 3.2. In Section 3.3, we specify the function $f$ for the maxima of sums and apply the theorems in Section 3.2 to derive a bound for the bootstrap LP pre-distance $\eta_n^*(\varepsilon, t)$. 

\subsection{Improved coherent Lindeberg interpolation method}
The standard Lindeberg interpolation \citep{lindeberg1922neue,chatterjee2005simple} bound the absolute difference of the expectation of $f(X_1, \ldots, X_n)$ and $f(Y_1, \ldots, Y_n)$ as
$\big| \mathbb{E} f(X_1, \ldots, X_n) - \mathbb{E} f(Y_1, \ldots, Y_n) \big| \le \sum_{i=1}^n \big| \mathbb{E} f(\mathbf{V}_{i-1}) - \mathbb{E} f(\mathbf{V}_{i}) \big|$, where $\mathbf{V}_i = (X_1, \ldots, X_i, Y_{i+1}, \ldots, Y_n)$.
However, as the analysis of the sum of interpolation difference on the right-hand side does not depend on the interpolation path $f(X_1, \ldots, X_n) = f(\mathbf{V}_n) \to f(\mathbf{V}_{n-1}) \to \cdots \to f(\mathbf{V}_1) \to f(\mathbf{V}_0) = f(Y_1, \ldots, Y_n)$ and this path is possibly the worst, the resulting upper bound may not be sharp. Motivated by this, the coherent Lindeberg interpolation method \citep{deng2017beyond} takes the average over all interpolation paths for permutation invariant $f$, that is, 
\begin{align*}
\big| \mathbb{E} f(X_1, \ldots, X_n) - \mathbb{E} f(Y_1, \ldots, Y_n) \big| \le \frac{1}{n!} \sum_{\sigma} \Big\{ \sum_{i=1}^n \big| \mathbb{E} f(\mathbf{V}_{\sigma, i-1}) - \mathbb{E} f(\mathbf{V}_{\sigma, i}) \big| \Big\},
\end{align*}
where $\mathbf{V}_{\sigma, i} = (X_{\sigma_1}, \ldots, X_{\sigma_i}, Y_{\sigma_{i+1}}, \ldots, Y_{\sigma_n})$. By taking average, the effect of bad interpolation paths may be reduced. After the interpolation argument, we apply Taylor series expansion to bound as $\mathbf{V}_{\sigma, i-1}$ and $\mathbf{V}_{\sigma, i}$ only differ at the $i$-th vector. In the Taylor series expansion, let the order of the derivative of $f$ in remainder term be $m^* > 2$. In Gaussian approximation considered in \cite{chernozhukov2013,chernozhukov2017central} where the first two moments are matched, $m^*=3$. If we consider higher order terms, $m^*$ should be greater than $3$. 

Let the $m$-th tensor valued derivative of $f$ with respect to the last input vector be 
$f^{(m)}(x_1, \ldots, x_n) := \big( \partial/\partial x_n \big)^{\otimes m} f(x_1, \ldots, x_n)$. In the sequel, we say $f^{(m)}$ is the $m$-th derivative of $f$. Denote $\boldsymbol{U}_{\sigma, i} = (X_{\sigma_1}, \ldots, X_{\sigma_{i-1}}, Y_{\sigma_{i + 1}}, \ldots, Y_{\sigma_n})$. Let $\mathbb{A}_{\sigma}$ be the average operator over all possible permutations $\sigma$, that is, $\mathbb{A}_{\sigma} h_{\sigma} = (n!)^{-1} \sum_{\sigma} h_{\sigma}$.

The full strength of the coherent Lindeberg interpolation method relies on the following Stability Condition \ref{cond:SC} imposed on $f^{(m^*)}(x_1, \ldots, x_n)$ and the Permutation Invariance Lemma \ref{lm:PI}. 

\begin{condition}[Stability Condition]\label{cond:SC}
There exists a function $\overline{f}^{(m^*)}(x_1, \ldots, x_n)$, a permutation invariant function $\overline{f}^{(m^*)}_{\max}(x_1, \ldots, x_n)$, and a nondecreasing function $g(\cdot)$ on $\mathbb{R}_+$ such that for a certain norm $\| \cdot \|$, it holds with probability $1$ that
\begin{align}
\label{1-cond:SC}\big| f^{(m^*)}(x_1, \ldots, x_{n-1}, t\xi ) \big| &\le g(\|t\xi\|) \overline{f}^{(m^*)}(x_1, \ldots, x_{n-1}, 0) \hbox{ and }
\\
\label{2-cond:SC} \overline{f}^{(m^*)}(x_1, \ldots, x_{n-1}, 0) 
& \le g(\|\xi\|) \overline{f}^{(m^*)}_{\max}(x_1, \ldots, x_{n-1}, \xi)
\end{align}
where $0 \le t \le 1$ and $\xi$ is either $X_i$ or $Y_i$ for any $1\le i \le n$.
\end{condition}

\begin{lemma}[Permutation Invariance Lemma]\label{lm:PI}
Let $\zeta_{i,k} = \delta_i X_k + (1 - \delta_i) Y_k$ with independent $\delta_i \sim \mathrm{Bernoulli}(\theta_{n,i})$. Define operator $\mathbb{A}_{\sigma, i} (\,\cdot\, ) := n^{-1} \sum_{i=1}^n \mathbb{A}_{\sigma} \mathbb{E} \big[\, \cdot \, \big| \boldsymbol{X}, \boldsymbol{Y}, \sigma, i\big]$. Then, for any permutation invariant function $f(x_1, \ldots, x_n)$, 
$ \mathbb{A}_{\sigma, i} I\{\sigma_i = k\} q_{n,i} f(\boldsymbol{U}_{\sigma,i}, \zeta_{i,\sigma_i})$
does not depend on $k$ for proper set $\{\theta_{n,i} \in [0,1], \forall i\}$ satisfying
\begin{align*}
(n-i)q_{n,i} \theta_{n,i} = i q_{n,i+1}(1- \theta_{n,i+1}) \qquad \forall \ 1\le i \le n-1.
\end{align*}
Consequently, for any function $g_k(\cdot, \cdot), 1\le k \le n$,
\begin{align}\label{1-lm:invariance}
& \mathbb{A}_{\sigma, i} q_{n,i}  \Big\langle f(\boldsymbol{U}_{\sigma,i}, \zeta_{i,\sigma_i}), \ g_{\sigma_i}( \boldsymbol{X} , \boldsymbol{Y}) \Big\rangle = \Big\langle \mathbb{A}_{\sigma, i} q_{n,i} f(\boldsymbol{U}_{\sigma,i}, \zeta_{i,\sigma_i}),\ \textstyle n^{-1} \sum_{k=1}^n g_k( \boldsymbol{X} , \boldsymbol{Y}) \Big\rangle.
\end{align}
\end{lemma}


\begin{remark}\label{rmk:PI}
This Lemma generalizes the Lemma 2 in \cite{deng2017beyond} in which $q_{n,i} \equiv 1$ so that $\theta_{n,i} = i/(n+1)$. In the proof of Theorem \ref{thm:DAT} in Section 3.3, we will encounter $q_{n,i} = (n+1-i)/(n+1)$, so $\theta_{n,i} = i/(n+2)$ should be used to apply this important lemma. 
\end{remark}

To understand the rationale of the coherent Lindeberg interpolation method, we shall look at a remainder term from the Taylor series expansion
\begin{align*}
\int_{0}^1 \frac{(1-\tau)^{m^*-1}}{(m^*-1)!} \bigg\{ \sum_{i=1}^n \mathbb{A}_{\sigma} \mathbb{E}  \Big\langle f^{(m^*)} (\boldsymbol{U}_{\sigma, i}, \tau X_{\sigma_i}), X_{\sigma_i}^{\otimes m^*} \Big\rangle \bigg\} d\tau.
\end{align*}
Suppose $\big\| f^{(m^*)} (x_1, \ldots, x_n) \big\|_1 \le D_n$ for any $(x_1, \ldots, x_n)$. If we apply H\"older's inequality directly, the main part is bounded by
\begin{align}\label{1-disc:CLIM}
\Big| \sum_{i=1}^n \mathbb{A}_{\sigma} \mathbb{E}  \Big\langle f^{(m^*)} ( \boldsymbol{U}_{\sigma, i}, \tau X_{\sigma_i}), X_{\sigma_i}^{\otimes m^*} \Big\rangle \Big| &\le \sum_{i=1}^n \mathbb{A}_{\sigma} \mathbb{E} \Big( \big\| f^{(m^*)} ( \boldsymbol{U}_{\sigma, i}, \tau X_{\sigma_i})\big\|_1 \cdot \big\| X_{\sigma_i}^{\otimes m^*} \big\|_{\max} \Big)
\\ \nonumber
&\le n D_n \cdot \frac{1}{n} \sum_{i=1}^n \mathbb{E} \max_j |X_{i,j} |^{m^*}.
\end{align}
The above bound is not ideal with the moment term $n^{-1} \sum_{i=1}^n \mathbb{E} \max_j |X_{i,j} |^{m^*}$ since the maximization of $X_{\sigma_i}^{\otimes}$ is taken before expectation $\mathbb{E}$. Our goal of using the coherent Lindeberg interpolation method is to reduce the moment term. With the help of the Stability Condition \ref{cond:SC} and Lemma \ref{lm:PI}, it follows that
\begin{align*}
&\Big| \sum_{i=1}^n \mathbb{A}_{\sigma} \mathbb{E}  \Big\langle f^{(m^*)} ( \boldsymbol{U}_{\sigma, i}, \tau X_{\sigma_i}), X_{\sigma_i}^{\otimes m^*} \Big\rangle \Big|
\\ \nonumber
&\le \sum_{i=1}^n \mathbb{A}_{\sigma} \mathbb{E}  \Big\langle \overline{f}^{(m^*)} ( \boldsymbol{U}_{\sigma, i}, 0), |X_{\sigma_i}|^{\otimes m^*}g(\|X_{\sigma_i} \|) \Big\rangle
\\ \nonumber
&= \sum_{i=1}^n \mathbb{A}_{\sigma} \Big\langle \mathbb{E}  \overline{f}^{(m^*)} (\boldsymbol{U}_{\sigma, i}, 0), \mathbb{E} |X_{\sigma_i}|^{\otimes m^*}g(\|X_{\sigma_i} \|) \Big\rangle
\\ \nonumber
&\le \sum_{i=1}^n \mathbb{A}_{\sigma} \Big\langle \mathbb{E}  \overline{f}^{(m^*)}_{\max} (\boldsymbol{U}_{\sigma, i}, \zeta_{i, \sigma_i}), \frac{\mathbb{E} |X_{\sigma_i}|^{\otimes m^*}g(\|X_{\sigma_i} \|) }{\mathbb{E} [ 1/g(\| \zeta_{i, \sigma_i} \|)] }  \Big\rangle
\\ \nonumber
&\le \sum_{i=1}^n \mathbb{A}_{\sigma} \Big\langle \mathbb{E}  \overline{f}^{(m^*)}_{\max} (\boldsymbol{U}_{\sigma, i}, \zeta_{i, \sigma_i}), \frac{1}{n}\sum_{k=1}^n \frac{\mathbb{E} |X_{k}|^{\otimes m^*}g(\|X_{k} \|) }{\mathbb{E} \big( [1/g(\| X_{k} \|)] \wedge [1/g(\| Y_{k} \|)] \big) }  \Big\rangle
\\ \nonumber
&\le \sum_{i=1}^n \mathbb{A}_{\sigma} \Big\| \mathbb{E}  \overline{f}^{(m^*)}_{\max} (\boldsymbol{U}_{\sigma, i}, \zeta_{i, \sigma_i}) \Big\|_1 \cdot \Big\| \frac{1}{n}\sum_{k=1}^n \frac{\mathbb{E} |X_{k}|^{\otimes m^*}g(\|X_{k} \|) }{\mathbb{E} \big( [1/g(\| X_{k} \|)] \wedge [1/g(\| Y_{k} \|)] \big) }  \Big\|_{\max}.
\end{align*}
In the above calculation, we (i) create independence between the two components in the inner product by \eqref{1-cond:SC} so that the expectation $\mathbb{E}$ can be taken first, (ii) introduce $\zeta_{i, \sigma_i}$ by \eqref{2-cond:SC} so that the Permutation Invariance Lemma \ref{lm:PI} is applicable, and (iii) apply H\"older's inequality to get the max norm at the end. For $g \le C$, we have the final bound
\begin{align*}
\textstyle n D_n \cdot C^2 \max_j \big( n^{-1} \sum_{i=1}^n \mathbb{E} |X_{i,j}|^{m^*} \big),
\end{align*}
which can be significantly smaller than the upper bound in \eqref{1-disc:CLIM}.

\subsection{General Comparison Theorems}\label{subsection:general_CT}

In this subsection, we present two general comparison theorems, the Weak Comparison Theorem \ref{thm:weak_comparison} and the Strong Comparison Theorem \ref{thm:strong_comparison}, to bound the weighted sum of interpolation differences
\begin{align}\label{def:Delta_nA}
\Delta_{n,\mathbb{A}} := \mathbb{A}_{\sigma} \sum_{i=1}^n q_{n, i} \Big( \mathbb{E} f(\boldsymbol{U}_{\sigma, i}, X_{\sigma_i}) - \mathbb{E} f( \boldsymbol{U}_{\sigma, i}, Y_{\sigma_i}) \Big),
\end{align}
for permutation invariant smooth function $f$. We note that when $q_{n,i} \equiv 1$ for all $i$, 
\begin{align*}
\Delta_{n,\mathbb{A}} = \mathbb{E} f(X_1, \ldots, X_n) - \mathbb{E} f(Y_1, \ldots, Y_n).
\end{align*}
In the end of this subsection, we will briefly compare the derived comparison theorems to other interpolation methods under $q_{n,i} \equiv 1$.

\begin{theorem}[Weak Comparison Theorem]\label{thm:weak_comparison}Let $(X_i, Y_i)\in \mathbb{R}^{p\times 2}$, $1\le i\le n$, be independent mean zero random matrices under expectation $\mathbb{E}$. Suppose $f$ is permutation invariant and its $m^*$-th derivative $f^{(m^*)}$ satisfies the  Stability Condition  \ref{cond:SC}. Suppose for $\{q_{n,i}, \forall i\}$ there exist $\{\theta_{n,i} \in [0,1], \forall i \}$ satisfying $(n- i)q_{n,i} \theta_{n,i} = i q_{n,i+1}(1- \theta_{n,i+1}), \ \forall \ 1\le i \le n-1$. Then, for $\Delta{n, \mathbb{A}}$ defined in \eqref{def:Delta_nA}
\begin{align}\label{1-thm:weak_comparison}
\Delta_{n,\mathbb{A}}  = \sum_{m=2}^{m^*-1} \frac{1}{m!} \sum_{i=1}^n q_{n,i} \mathbb{A}_{\sigma} \Big\langle \mathbb{E} f^{(m)} ( \boldsymbol{U}_{\sigma, i}, 0), \mathbb{E} X_{\sigma_i}^{\otimes m} - \mathbb{E} Y_{\sigma_i}^{\otimes m} \Big\rangle + \mathrm{Rem}_1,
\end{align}
where $|\mathrm{Rem}_1|$ is bounded by
\begin{align}\label{2-thm:weak_comparison}
|\mathrm{Rem}_1| \le
 \bigg\langle & \sum_{i=1}^n q_{n,i} \mathbb{A}_{\sigma} \, \mathbb{E}  \overline{f}^{(m^*)}_{\max} ( \boldsymbol{U}_{\sigma, i}, \zeta_{i, \sigma_i}), \frac{1}{n}\sum_{k=1}^n \frac{\mathbb{E} |X_k|^{\otimes m^*} g(\|X_k\|) + \mathbb{E} |Y_k|^{\otimes m^*} g(\|Y_k\|)}{(m^*)! \, \mathbb{E} [1/g(\|X_k\|)] \wedge \mathbb{E} [1/g(\|Y_k\|)]} \bigg\rangle.
\end{align}
Here $\zeta_{i, \sigma_i} = \delta_i X_{\sigma_i} + (1- \delta_i) Y_{\sigma_i}$ where $\{ \delta_i \overset{ind.}{\sim} \mathrm{Bernoulli}(\theta_{n, i}), \forall i \}$ are independent of $\{(X_i, Y_i), \forall i\}$.
\end{theorem}



If we further assume the permutation invariance of $f^{(m)}$ for $2 \le m < m^*$, we have the following Strong Comparison Theorem. However, it \textit{only} works for $m^*=3$ and $4$.

\begin{theorem}[Strong Comparison Theorem]\label{thm:strong_comparison} Let $(X_i, Y_i)\in \mathbb{R}^{p\times 2}$, $1\le i\le n$, be independent mean zero random matrices under expectation $\mathbb{E}$. Let $m^*=3$ or $4$. Suppose $f$ and its derivatives $f^{(m)}$ of order $m = 2, \ldots, m^*-1$ are permutation invariant and the $m^*$-th derivative $f^{(m^*)}$ satisfies the Stability Condition \ref{cond:SC}. Suppose for $\{q_{n,i}, \forall i\}$ there exist $\{\theta_{n,i} \in [0,1], \forall i \}$ satisfying $(n- i)q_{n,i} \theta_{n,i} = i q_{n,i+1}(1- \theta_{n,i+1}), \ \forall \ 1\le i \le n-1$. Then, for $\Delta_{n,\mathbb{A}}$ defined in \eqref{def:Delta_nA},
\begin{align}\label{1-prop:strong_comparison}
\Delta_{n,\mathbb{A}} &= \sum_{m=2}^{m^*-1} \frac{1}{m!}\Big\langle \sum_{i=1}^n q_{n,i} \mathbb{A}_{\sigma} \mathbb{E} f^{(m)} ( \boldsymbol{U}_{\sigma, i}, \zeta_{i, \sigma_i}), \frac{1}{n} \sum_{k=1}^n \big(\mathbb{E} X_k^{\otimes m} - \mathbb{E} Y_k^{\otimes m} \big) \Big\rangle + \mathrm{Rem}_2,
\end{align}
where, with $\zeta_{i, \sigma_i}$ as in Theorem \ref{thm:weak_comparison}, $|\mathrm{Rem_2}|$ is bounded by
\begin{align}\label{2-prop:strong_comparison}
|\mathrm{Rem}_2| 
&\le \frac{2^{m^*} - m^*-1}{(m^*)!} \bigg\langle \sum_{i=1}^n q_{n,i} \mathbb{A}_{\sigma} \mathbb{E}  \overline{f}^{(m^*)}_{\max} ( \boldsymbol{U}_{\sigma, i}, \zeta_{i, \sigma_i}), \ \mu_g^{(m^*)}  \bigg\rangle.
\end{align}
Here with $G_k = \big( \mathbb{E} [1/ g(\|X_k\|)] \big) \wedge \big( \mathbb{E} [1/ g(\|Y_k\|)] \big)$
\begin{align*}
\mu_g^{(m^*)} &:= \bigg[ \Big(\frac{1}{n} \sum_{k=1}^n \frac{\mathbb{E} |X_k|^{m^*} g(\|X_k\|)}{G_k} \Big)^{1/m^*} \bigg]^{\otimes m^*} + \bigg[ \Big(\frac{1}{n} \sum_{k=1}^n \frac{\mathbb{E} |X_k|^{m^*} g(\|Y_k\|)}{G_k} \Big)^{1/m^*} \bigg]^{\otimes m^*}
\cr
&\quad + \bigg[ \Big(\frac{1}{n} \sum_{k=1}^n \frac{\mathbb{E} |Y_k|^{m^*} g(\|X_k\|)}{G_k} \Big)^{1/m^*} \bigg]^{\otimes m^*} + \bigg[ \Big(\frac{1}{n} \sum_{k=1}^n \frac{\mathbb{E} |Y_k|^{m^*} g(\|Y_k\|)}{G_k} \Big)^{1/m^*} \bigg]^{\otimes m^*}
\end{align*}
\end{theorem}


While Theorem \ref{thm:weak_comparison} is new, Theorem \ref{thm:strong_comparison} is a generalized version of Theorem 4 in \cite{deng2017beyond} in which $q_{n,i}$ is set to be fixed constant $1$, so that $\Delta_{n, \mathbb{A}} = \mathbb{E}f(X_1, \ldots, X_n) - \mathbb{E} f(Y_1, \ldots, Y_n)$. Recall from Remark \ref{rmk:PI}, $\zeta_{i, \sigma_i} = \delta_i X_{\sigma_i} + (1 - \delta_i) Y_{\sigma_i}$ with $\delta_i \overset{ind.}{\sim} \mathrm{Bernoulli}(i/(n+1))$ for $q_{n,i} \equiv 1$.

\medskip
\textit{Comparison with other interpolation methods.}
We compare the improved coherent Lindeberg interpolation method to some common interpolation methods. In common interpolation methods, the quantity of interest is exclusively $\Delta_{n, \mathbb{A}} = \mathbb{E}f(X_1, \ldots, X_n) - \mathbb{E} f(Y_1, \ldots, Y_n)$, that is, $q_{n, i} \equiv 1$. We focus on this case in the discussion.

By Taylor series expansion, we know the standard Lindeberg interpolation \citep{lindeberg1922neue,chatterjee2005simple} is only able to yield a similar result to the Weak Comparison Theorem \ref{thm:weak_comparison} with operator $\mathbb{A}_{\sigma}$ removed and the remainder term $\mathrm{Rem}^{\mathrm{(std)}}$ bounded by
\begin{align*}
|\mathrm{Rem}^{\mathrm{(std)}}| & \lesssim \frac{1}{(m^*)!} \sum_{i=1}^n \Big\langle \mathbb{E}  \overline{f}^{(m^*)} ( \boldsymbol{U}_{\sigma, i}, 0), \mathbb{E} |X_i|^{\otimes m^*} + \mathbb{E} |Y_i|^{\otimes m^*} \Big\rangle
\cr
& \lesssim \max_i \big\| \mathbb{E}  \overline{f}^{(m^*)} ( \boldsymbol{U}_{\sigma, i}, 0)\big\|_1 \cdot \Big[\frac{1}{n}\sum_{k=1}^n \max_j \big( \mathbb{E} |X_{k,j}|^{m^*} + \mathbb{E} |X_{k,j}|^{m^*} \big) \Big].
\end{align*}
It is clearly sub-optimal if $f^{(m^*)}$ satisfies the Stability Condition \ref{cond:SC} with $g$ properly controlled, in view of \eqref{2-thm:weak_comparison} in Theorem \ref{thm:weak_comparison} and \eqref{2-prop:strong_comparison} in Theorem \ref{thm:strong_comparison} where the maximizations over $j$ are both taken after the average over $k$. 

The Slepian's `smart' interpolation \citep{chernozhukov2013,chernozhukov2015comparison}, resolves this issue and yields a result similar to the Weak Comparison Theorem \ref{thm:weak_comparison} but only for $m^*=3$. Instead of taking the interpolation paths as in the coherent Lindeberg interpolation, the Slepian's `smart' interpolation takes a continuous path $f \big(Z_1(\theta), \ldots, Z_n(\theta) \big)$ with $\theta $ goes from $0$ to $\pi/2$ where $Z_i = X_i \cos \theta  + Y_i \sin \theta$. It implies
\begin{align*}
& \big| \mathbb{E} f(X_1, \ldots, X_n) - \mathbb{E} f(Y_1, \ldots, Y_n) \big| 
\cr
& = \big| \mathbb{E} \int_{0}^{\pi/2} \mathrm{d} \, f(Z_1, \ldots, Z_n) \big|
\cr
& = \sum_{i=1}^n \mathbb{E} \int_{0}^{\pi/2} \Big\langle \Big( \frac{\partial}{\partial Z_i} \Big) f(Z_1, \ldots, Z_n), - X_i \sin \theta + Y_i \cos \theta \Big\rangle \mathrm{d}\, \theta 
\cr
& = \sum_{i=1}^n \int_{0}^{\pi/2} \Big\langle \mathbb{E} \Big(\frac{\partial}{\partial Z_i}\Big)^{\otimes 2} f(Z_1, \ldots, Z_n) \big|_{Z_i = 0}, \ \mathbb{E} Z_i \otimes \big( - X_i \sin \theta + Y_i \cos \theta \big) \Big\rangle \mathrm{d}\, \theta + \mathrm{Rem}.
\end{align*}
This yields the second-order moment comparison term $\mathbb{E} \big[ Z_i \otimes ( - X_i \sin \theta + Y_i \cos \theta) \big] = \cos \theta \sin \theta \big[ \mathbb{E} Y_i^{\otimes 2} - \mathbb{E} X_i^{\otimes 2} \big]$. However, such comparison is only possible for $m^*=3$. If we have the third moment match $\mathbb{E} X_i^{\otimes 3} = \mathbb{E} Y_i^{\otimes 3}$, a higher order Taylor series expansion produces $\mathbb{E} \big[ Z_i^{\otimes 2} \otimes (-X_i \sin \theta + Y_i \cos \theta) \big] = \sin^2 \theta \cos \theta \mathbb{E} Y_i^{\otimes 3} - \cos^2 \theta \sin \theta \mathbb{E} X_i^{\otimes 3}$; this term fails to vanish for equal third moments. In contrast, the Weak Comparison Theorem \ref{thm:weak_comparison} can be used to compare arbitrarily many moments as it holds for any $m^* \ge 3$.

Moreover, we point out that the Weak Comparison Theorem \ref{thm:weak_comparison}, the classical Lindeberg interpolation and the Slepian's `smart' interpolation can only manage moment comparison at individual level, that is, they compare $\mathbb{E} X_i^{\otimes m}$ and $\mathbb{E} Y_i^{\otimes m}$ for each $i$. However, the Strong Comparison Theorem \ref{thm:strong_comparison} indicates that this is not necessary at least for $m^*=3$ and $4$  --- it suffices to compare their average moments, $n^{-1} \sum_{i=1}^n \mathbb{E} X_i^{\otimes m}$ and $n^{-1} \sum_{i=1}^n \mathbb{E} Y_i^{\otimes m}$. This is a helpful feature to deal with heterogeneous $\mathbb{E} X_i^{\otimes m} - \mathbb{E} Y_i^{\otimes}$, e.g., in the case of empirical bootstrap with non-i.i.d. $\{X_i\}$.

Lastly, it is worth mentioning that \cite{chernozhukov2019improved} also modifies the original coherent Lindeberg interpolation in \cite{deng2017beyond} and proposed the so-called iterative randomized Lindeberg interpolation method to sharpen their consistency rates for bootstrap under sub-Gaussian condition \ref{ex-gaussian}, discussed in Example 1. However, their method appears to require some strong conditions, e.g. an $8$-th moment condition on data $\{X_i\}$, and is not as general as the improved coherent Lindeberg interpolation method in this paper.

\subsection{Comparison Theorem for the maxima of independent sums}\label{subsection:maxima}
We apply the improved coherent Lindeberg interpolation method to study the LP pre-distance between the maxima of sums $T_n$ and its counterpart $T_n^Y = \max_j \sum_{i=1}^n Y_i/\sqrt{n}$,
\begin{align}\label{def:general_LP}
\eta_n(\varepsilon, t; T_n, T_n^Y) &= \eta_n^{(\mathbb{P})}(\varepsilon, t; T_n, T_n^Y) 
\\ \nonumber 
&:= \max \Big[0, \ \mathbb{P} \big\{ T_n \le t -\varepsilon \big\} - \mathbb{P} \big\{ T_n^Y \le t \big\}, \mathbb{P} \big\{ T_n^{Y} \le t -\varepsilon \big\} - \mathbb{P} \big\{ T_n  \le t \big\} \Big].
\end{align}
It is connected to the Kolmogorov-Smirnov distance between $T_n$ and $T_n^Y$ via 
\begin{align}\label{KS-bound}
\sup_t \big| \mathbb{P}\{T_n \le t\} - \mathbb{P} \{T_n^Y \le t\} \big| \le \sup_t \eta_n(\varepsilon, t; T_n, T_n^Y) + \omega_n(\varepsilon; T_n) \wedge \omega_n(\varepsilon; T_n^Y) \quad \forall \varepsilon >0,
\end{align}
where the anti-concentration $\omega_n(\varepsilon; \xi)$ of random variable $\xi$ is defined as
\begin{align}\label{def:anti-concentration}
\omega_n(\varepsilon; \xi) = \sup_{t} \mathbb{P}\{t -  \varepsilon \le \xi < t \}.
\end{align}

For the bootstrap LP pre-distance $\eta_n^*(\varepsilon, t)$ in \eqref{def:bootstrap_LP}, we let $\{X_i^0\}$ be an independent copy of $\{X_i\}$, and define $T_n^0 = \max_j \sum_{i=1}^n X_i^0/\sqrt{n}$. Let the probability measure in \eqref{def:general_LP} be the bootstrap probability measure $\mathbb{P}^* = \mathbb{P} \{ \, \cdot \, |X_i, \forall i\}$, then we have
$\eta_n^{(\mathbb{P^*})}(\varepsilon, t; T^0_n, T_n^*)  = \eta_n^*(\varepsilon, t)$.
Therefore, it suffices to focus our analysis on $\eta_n(\varepsilon, t; T_n, T_n^Y) = \eta_n^{(\mathbb{P})}(\varepsilon, t; T_n, T_n^Y)$ for general $\mathbb{P}$.

To apply the improved coherent interpolation method, we shall specify the smooth function $f$. Indeed, for a smooth function $h(\cdot)$ with $h(t)=1$ for $t \le 0$ and $h(t) = 0$ for $t \ge 1$, we may use $h\big(2(T_n - t)/\varepsilon\big)$ to approximate the indicator function $I\{T_n \le t\}$ and they only differ at a small interval $\{t < T_n < t + \varepsilon/2 \}$ with length $\varepsilon/2$. Moreover, since the maximum $T_n$ is not a smooth function of $(X_1, \ldots, X_n)$, we approximate it with the softmax function 
$F_\beta(z) = \beta^{-1} \log\big( e^{\beta z_1}+\cdots+e^{\beta z_p}\big)$ with $z=(z_1,\ldots,z_p)^\top$. As a result, the LP pre-distance $\eta_n(\varepsilon, t; T_n, T_n^Y)$ can be bounded by
\begin{align*}
\eta_n(\varepsilon, t; T_n, T_n^Y) \le  \big| \mathbb{E} f(X_1, \ldots, X_n) - \mathbb{E} f(Y_1, \ldots, Y_n) \big|,
\end{align*}
where with $z = \sum_{i=1}^n x_i/\sqrt{n}$
\begin{align}\label{def:f}
f(x_1, \ldots, x_n) = f_t(x_1, \ldots, x_n) := h \big( 2 \varepsilon^{-1} F_{\beta}(z) - 2t/\varepsilon + 1 \big)
\end{align}
and $\beta \ge 2(\log p)/\varepsilon$. Two key properties of this function $f$, studied in \cite{deng2017beyond}, are given in Proposition \ref{prop:H} below.

\begin{proposition}[\cite{deng2017beyond}]\label{prop:H} Let $h_0(\cdot)$ be a smooth function and $z = \sum_{i=1}^n x_i/\sqrt{n}$.  
Then, there exist functions $H^{(m)}_{\varepsilon, \beta}(z)$ for $m \ge 1$ that satisfy
\begin{align}\label{f_stability-1}
& \big|n^{m/2}(\partial/\partial x_n)^{\otimes m} h_0( 2 \varepsilon^{-1} F_{\beta}(z) ) \big| \le H^{(m)}_{\varepsilon, \beta}(z), 
\\ \nonumber
& \big\|H^{(m)}_{\varepsilon, \beta}(z) \big\|_1 \le C_{h,m} \max \big\{\varepsilon^{-m}, \varepsilon^{-1} \beta^{m-1} \big\},
\end{align} 
and
\begin{align}\label{f_stability-2}
e^{-2m\|t\|_{\infty}\beta}H^{(m)}_{\varepsilon,\beta}(z+t) \le H^{(m)}_{\varepsilon,\beta}(z)
\le e^{2m\|t\|_{\infty}\beta}H^{(m)}_{\varepsilon,\beta}(z+t).
\end{align}
\end{proposition} 
This proposition essentially allows us to find function triplet $(\overline{f}^{(m^*)}, \overline{f}^{(m^*)}_{\max}, g)$ to meet the Stability Condition \ref{cond:SC} and apply the comparison theorems in Section \ref{subsection:general_CT}. For $f = f_t$ in \eqref{def:f}, \cite{deng2017beyond} considered $\overline{f}^{(m^*)}(x_1, \ldots, x_n) = n^{-m^*/2} H^{(m^*)}_{\varepsilon,\beta}(z)$, $\overline{f}^{(m^*)}_{\max}(x_1, \ldots, x_n) = n^{-m^*/2} H^{(m^*)}_{\varepsilon,\beta} (z)$ and $g(\|\xi\|_{\infty}) = \exp \big(2m^* \max_j |\xi_j| \beta/\sqrt{n} \big)$. However, by the definition of $f$ in \eqref{def:f} and the fact that $h$ is constant except on $[0,1]$, we observe $f^{(m)}(x_1, \ldots, x_n) = f^{(m)}(x_1, \ldots, x_n) I\{t - \varepsilon \le \max_j z_j \le t\}$ for $m \ge 1$. As a result, for $\max_j |\xi| \le c_n \sqrt{n}/\beta$, that is, $\max_{i,j} |X_{i,j}| \vee |Y_{i,j}| \le c_n \sqrt{n}/\beta$, it is better to consider
\begin{align}\label{f_stability_bounded}
\overline{f}^{(m^*)}(x_1, \ldots, x_n) &:= n^{-m^*/2} H^{(m^*)}_{\varepsilon,\beta}(z) I\big\{t - \varepsilon - \frac{c_n}{\beta} \le \max_j z_j  \le t + \frac{c_n}{\beta} \big\},
\\ \nonumber
\overline{f}^{(m^*)}_{\max}(x_1, \ldots, x_n) &:= n^{-m^*/2} H^{(m^*)}_{\varepsilon,\beta} (z) I\big\{t - \varepsilon - \frac{2c_n}{\beta} \le \max_j z_j \le t + \frac{2c_n}{\beta} \big\},
\\ \nonumber
g(\|\xi\|_{\infty}) &:= \exp \big(2m^* \max_j |\xi_j| \beta/\sqrt{n} \big).
\end{align}
It is easy to show such $(\overline{f}^{(m^*)}, \overline{f}^{(m^*)}_{\max}, g )$ satisfies the Stability Condition \ref{cond:SC}. As a function of the sum $\sum_{i=1}^n x_i$, $f^{(m)}(x_1, \ldots, x_n)$ must be permutation invariant for all $m \ge 1$. Overall, this implies that both the Weak Comparison Theorem \ref{thm:weak_comparison} and the Strong Comparison Theorem \ref{thm:strong_comparison} are applicable to the $f = f_t$ defined in \eqref{def:f}.

The new configuration of $(\overline{f}^{(m)}, \overline{f}^{(m^*)}_{\max})$ in \eqref{f_stability_bounded} plays a crucial role in contributing to the improvement of our final results over \cite{deng2017beyond} and other works. The sharp consistency rates in Theorems \ref{thm:main}, \ref{thm:main-full} and \ref{thm:main-exact} can all be traced back to the extra possibly small indicator functions in \eqref{f_stability_bounded}.

With Stability Condition \ref{cond:SC} satisfied via \eqref{f_stability_bounded}, we will repeatedly apply the Strong Comparison Theorem \ref{thm:strong_comparison} to bound the LP pre-distance $\eta_n(\varepsilon, t; T_n, T_n^Y)$ in \eqref{def:general_LP}. This is stated as the following distributional approximation theorem.

\begin{theorem}[Distributional Approximation Theorem]\label{thm:DAT}
Let $(X_i, Y_i)\in \mathbb{R}^{p\times 2}$, $1\le i\le n$, be independent mean zero random matrices under expectation $\mathbb{E}$. Let $M_{m^*}$ be as in \eqref{average_moments} and $M_{m^*, Y}$ the counterpart of $M_{m^*}$ for $\{Y_i\}$. Consider $m^* = 3$ or $4$.
Suppose there exist independent random variables $\{\mathfrak{W}_i, \forall i\}$, also independent of $\{(X_i, Y_i), \forall i\}$, such that $\mathbb{P}\big\{ \max_j |X_{i,j}| \vee |Y_{i,j}| \le |\mathfrak{W}_i| \sqrt{n} \varepsilon/(\log (np))$ and $|\mathfrak{W}_i| \le c \log(np) \ \forall i \big\} = 1$ and $\max_i \Big( \mathbb{E} \exp\{4m^* |\mathfrak{W}_i| \} \big/ \mathbb{E} \exp\{-4m^* |\mathfrak{W}_i| \} \Big) \le C_{c, m^*}^{\mathfrak{W}}$ for a constant $C_{c_0, m^*}^{\mathfrak{W}}$. Let $\omega_n(\varepsilon; \xi)$ be as in \eqref{def:anti-concentration} and define 
\begin{align*}
K_{n, m^*}(\varepsilon) &= \displaystyle \sum_{m=2}^{m^*-1} C_{c, m} \frac{(\log (np) )^{m-1}}{n^{m/2-1} \varepsilon^m} \Big\| \frac{1}{n}\sum_{i=1}^n \mathbb{E} X_i^{\otimes m} - \frac{1}{n}\sum_{i=1}^n \mathbb{E} Y_i^{\otimes m} \Big\|_{\max}
\\ \nonumber
& \qquad  + C_{c, m^*}^{\mathfrak{W}} C_{m^*} \frac{(\log (np))^{m^*-1}}{ n^{m^*/2-1}\varepsilon^{m^*}} (M_{m^*}^{m^*} + M_{m^*, Y}^{m^*}).
\end{align*}
Then, the LP pre-distance $\eta_n(\varepsilon, t; T_n, T_n^Y)$ is bounded by 
\begin{align}\label{1-thm:DAT}
\eta_n(\varepsilon, t; T_n, T_n^Y) &= \max \Big[0, \ \mathbb{P} \big\{ T_n \le t -\varepsilon \big\} - \mathbb{P} \big\{ T_n^Y \le t \big\}, \mathbb{P} \big\{ T_n^{Y} \le t -\varepsilon \big\} - \mathbb{P} \big\{ T_n  \le t \big\} \Big]
\\ \nonumber
&\le  K_{n, m^*}(\varepsilon) \min \bigg\{1, \frac{\omega_n \big( (4c+3) \varepsilon; T_n \big) \wedge \omega_n\big( (4c+3) \varepsilon; T_n^Y \big)}{ \big[ 1-K_{n, m^*}(\varepsilon) \big]_+ } \bigg\}
\end{align}
\end{theorem}


Suppose we do not take into account the indicator functions in \eqref{f_stability_bounded} and use the same configuration of $(\overline{f}^{(m^*)}, \overline{f}^{(m^*)}_{\max}, g)$ as in \cite{deng2017beyond}, it is only possible to show
\begin{align*}
&\eta_n(\varepsilon, t; T_n, T_n^Y) \le  K_{n, m^*}(\varepsilon).
\end{align*}
The improvement here in Theorem \ref{thm:DAT} is the extra minimization factor on the right-hand side of \eqref{1-thm:DAT}, which can be significantly smaller than $1$ for small $\varepsilon$. 

In Theorem \ref{thm:DAT}, random variables $\mathfrak{W}_i$'s are introduced to essentially allow, e.g., $Y_i = W_i X_i$. This is very helpful when we derive theoretical results for the multiplier bootstrap. In empirical bootstrap $\mathfrak{W}_i$ can be simply set to be a fixed constant. For unbounded multipliers and $\{(X_i, Y_i)\}$, we may employ truncation arguments and the pseudo-triangle inequality of LP pre-distance 
\begin{align}\label{triangle_inequality}
\eta_n(\varepsilon, t; T_n, T_n^Y) &\le \sup_t \eta_n(\varepsilon_1, t; T_n, \widetilde{T}_n) + \sup_t \eta_n(\varepsilon_2, t; \widetilde{T}_n, \widetilde{T}_n^Y) + \sup_t \eta_n(\varepsilon_3, t; \widetilde{T}_n^Y, T_n^Y),
\end{align}
where $\varepsilon_1 + \varepsilon_2 + \varepsilon_3 \le \varepsilon$ and $\widetilde{T}_n$ and $\widetilde{T}_n^Y$ are the truncated versions of $T_n$ and $T_n^Y$ respectively. The truncation effects, $\sup_t \eta_n(\varepsilon_1, t; T_n, \widetilde{T}_n)$ and $\sup_t \eta_n(\varepsilon_3, t; \widetilde{T}_n^Y, T_n^Y)$, can be controlled by certain tail probabilities of $\{X_i\}$ and $\{Y_i\}$. 

It remains to bound the anti-concentration of $T_n$ or $T_n^Y$ to make Theorem \ref{thm:DAT} immediately applicable as $\omega_n ( \varepsilon; T_n ) \wedge \omega_n(\varepsilon; T_n^Y )$ appears in \eqref{1-thm:DAT}. The following theorem takes care of it by giving an anti-concentration bound for general $T_n$.

\begin{theorem}\label{thm:gAC}
Let $X_i\in \mathbb{R}^{p}$, $1\le i\le n$, be independent mean zero random vectors under expectation $\mathbb{E}$. Let $\overline{\sigma}$ be as in \eqref{def:soft_minimum}. Let $\widetilde{X}_{i,j} = X_{i,j}I\{|X_{i,j}| \le a_n\} - \mathbb{E} X_{i,j}I\{|X_{i,j}| \le a_n\}$ and $a_n = c_0 \sqrt{n}\varepsilon/(\log (np))$ for a fixed constant $c_0$. Then, the anti-concentration $\omega_n(\varepsilon; T_n)$, defined in \eqref{def:anti-concentration}, is bounded by
\begin{align}\label{1-thm:gAC}
\omega_n(\varepsilon; T_n) &\le C_{c_0} \bigg\{ \frac{(\log (np))^3(\log p)^{1/2}}{n} \frac{M_4^4}{\varepsilon^3 \overline{\sigma}} + \frac{\varepsilon}{\overline{\sigma}}\sqrt{\log p} \bigg\} + 2\mathbb{P} \Big[ \Big\| \frac{1}{\sqrt{n}}\sum_{i=1}^n \big( X_{i,j} - \widetilde{X}_{i, j} \big) \Big\|_{\infty} \ge \frac{\varepsilon}{2} \Big].
\end{align}	
Moreover, it holds when $ \rho_n := \displaystyle \frac{4(\log (np))^3}{c_0^3 \cdot n \, \varepsilon^4} \mathbb{E} \max_j \frac{1}{n}\sum_{i=1}^n X_{i,j}^4I\{|X_{i,j}| \ge a_n\} \le 1$ that
\begin{align}\label{2-thm:gAC}
\mathbb{P} \bigg[ \Big\| \frac{1}{\sqrt{n}}\sum_{i=1}^n \big( X_{i,j} - \widetilde{X}_{i, j} \big) \Big\|_{\infty} \ge \frac{\varepsilon}{2} \bigg] &\le 
\mathbb{P} \bigg[ \Big\| \frac{1}{\sqrt{n}}\sum_{i=1}^n X_{i,j}I\{|X_{i,j}| \ge a_n\} \Big\|_{\infty} \ge \frac{\varepsilon}{4} \bigg]  \le \rho_n.
\end{align}
\end{theorem}


The Gaussian anti-concentration $\omega_n(\varepsilon; \max_{1 \le j \le p} \xi_j)$ with $\xi$ being a Gaussian vector has been studied in many works, e.g., \cite{nazarov2003maximal,klivans2008learning,chernozhukov2015comparison}. The sharpest result for general $\xi$ is probably due to \cite{deng2017beyond}. Let $\xi_j \sim N(\mu_j, \sigma_j^2)$ and $\overline{\sigma}$ be defined as in \eqref{def:soft_minimum}, Theorem 10 of \cite{deng2017beyond} proved that
\begin{align*}
\omega_n \big(\varepsilon; \max_{1 \le j \le p} \xi_j \big) \le \frac{\varepsilon}{\overline{\sigma}} ( 2 + \sqrt{2 \log p})
\end{align*}
and it is rate optimal. The proof strategy for Theorem \ref{thm:gAC} is then to take advantage of this Gaussian anti-concentration bound by approximating $T_n$ with $T_n^Y$ where $Y_i = W_i X_i$ and $W_i$ has a Gaussian component.

\section{Simulation Results}\label{section:simulation}
In this section, we study the numerical performance of conservative bootstrap for finite samples. The following numerical study of the conservative bootstrap procedure serves as a complementary material to support our theoretical results in Section 2, that is, the conservative bootstrap, using either the empirical bootstrap (EB) or the Mammen's multiplier bootstrap (MB), should perform well in yielding small one-sided coverage error $\eta_{n, \alpha}^*$ defined in \eqref{conservative-bootstrap-error}. Although our theoretical results do not cover the Gaussian multiplier bootstrap (GB) or the Rademacher multiplier bootstrap (RB) in general, we include them to compare their numerical performance with MB and EB.

We generate $n$ independent $p$-dimensional vectors $\{X_1, \ldots, X_n\}$ in a Gaussian copula model: First draw $n$ i.i.d.~Gaussian vectors $\{Y_1, \ldots, Y_n\}$ from $\mathcal{N}(0, \Sigma)$ with marginal distributions $\mathcal{N}(0,1)$, and then let $X_{i,j}$ be such that $F(X_{i,j}) = \Phi(Y_{i,j})$ for all $i,j$ where $\Phi$ is the cdf of $\mathcal{N}(0,1)$ and $F$ is the cdf of gamma distribution with unit scale and shape parameter $1$. We set sample size $n=200$ and dimension $p=10^3$. Four different settings on covariance matrix $\Sigma$ are considered: (a) $\Sigma_{j,k} = I\{ j = k \}$, (b) $\Sigma_{j,k} = 0.2^{|j-k|}$, (c) $\Sigma_{j,k} = 0.8^{|j-k|}$ and (d) $\Sigma_{j,k} = 0.8 + 0.2 I\{j=k\}$. Note that in these settings the Gaussian multiplier bootstrap and the Rademacher multiplier bootstrap do not match the third moment as in \eqref{third_moment_match}. To compare the one-sided coverage error $\eta_{n, \alpha}^*$ defined in \eqref{conservative-bootstrap-error}, we compute the coverage probabilities $\mathbb{P}\{T_n \le 1.01 t_{\alpha}^*\}$ for these bootstrap schemes as follows: (i) In each experiment, generate $K = 10^4$ sets of data $\{X_1^{(k)}, \ldots, X_n^{(k)}\}$ where $1 \le k \le K$; (ii) for each set of data $\{X_1^{(k)}, \ldots, X_n^{(k)}\}$, compute the maximum $T_n$, denoted by $T_n^{(k)}$, and apply the aforementioned four bootstrap schemes (GB, MB, RB, EB) with $B = 10^3$ bootstrap samples to obtain the corresponding bootstrap quantiles $t_{\alpha}^*$, denoted by $(t_{\alpha}^*)^{(k)}$; (iii) find the relative frequency of $T_n^{(k)} \le 1.01(t_{\alpha}^*)^{(k)}$, that is, we approximate the conservative coverage probability $\mathbb{P}\{ T_n \le 1.01 t_{\alpha}^*\}$ with $K^{-1} \sum_{k=1}^K I\{ T_n^{(k)} \le 1.01 (t_{\alpha}^*)^{(k)}\}$.
We also give the relative frequency of $T_n^{(k)} \le (t_{\alpha}^*)^{(k)}$ to approximate the exact coverage probability $\mathbb{P}\{ T_n \le t_{\alpha}^*\}$. We consider $\alpha = 0.05$, so the targeted coverage probability is $0.95$. The simulation results are presented in Table \ref{tab:relative_frequency}.

\begin{table}[htb]
\centering
	\begin{tabular}{c||c c c c|c c c c}
		\hline 
		 & \multicolumn{4}{c}{$\mathbb{P}\{ T_n \le 1.01t_{\alpha}^*\}$ } & \multicolumn{4}{c}{$\mathbb{P}\{ T_n \le t_{\alpha}^*\}$} \\
		\hline
		Exp. & GB & MB & RB & EB & GB & MB & RB & EB \\
		\hline
		(a) & 0.9334 & 0.9640 & 0.8853 & 0.9818 & 0.9226 & 0.9567 & 0.8724 & 0.9780 \\
		\hline
		(b) & 0.9330 & 0.9631 & 0.8866 & 0.9810 & 0.9219 & 0.9569 & 0.8709 & 0.9782 \\
		\hline
		(c) & 0.9325 & 0.9616 & 0.8934 & 0.9768 & 0.9242 & 0.9562 & 0.8819 & 0.9732 \\
		\hline
		(d) & 0.9434 & 0.9707 & 0.9145 & 0.9826 & 0.9357 & 0.9664 & 0.9058 & 0.9798 \\
		\hline
	\end{tabular}
	\vspace{0.5ex}
	\caption{Simulated relative frequencies of $\{ T_n \le 1.01t_{\alpha}^* \}$ (conservative) and $\{ T_n \le t_{\alpha}^* \}$ (exact).}
	\label{tab:relative_frequency}
\end{table}

From Table \ref{tab:relative_frequency}, we can see it clearly that the conservative one-sided coverage error $\eta_{n, \alpha}^*$ diminishes in all settings as the simulated relative frequencies are greater than $95\%$. The two-sided coverage errors are in general small (less than $3\%$ for GB, MB, EB) except for the Rademacher multiplier bootstrap. For more detailed comparison of the exact bootstrap procedures, we refer readers to \cite{deng2017beyond} and \cite{chernozhukov2019improved}.

\begin{figure}[hbt]
	\centering
	\subfigure[$\Sigma_{j,k} = I\{ j = k \}$]{
		\label{fig:exp1} 
		\includegraphics[width=0.45\textwidth]{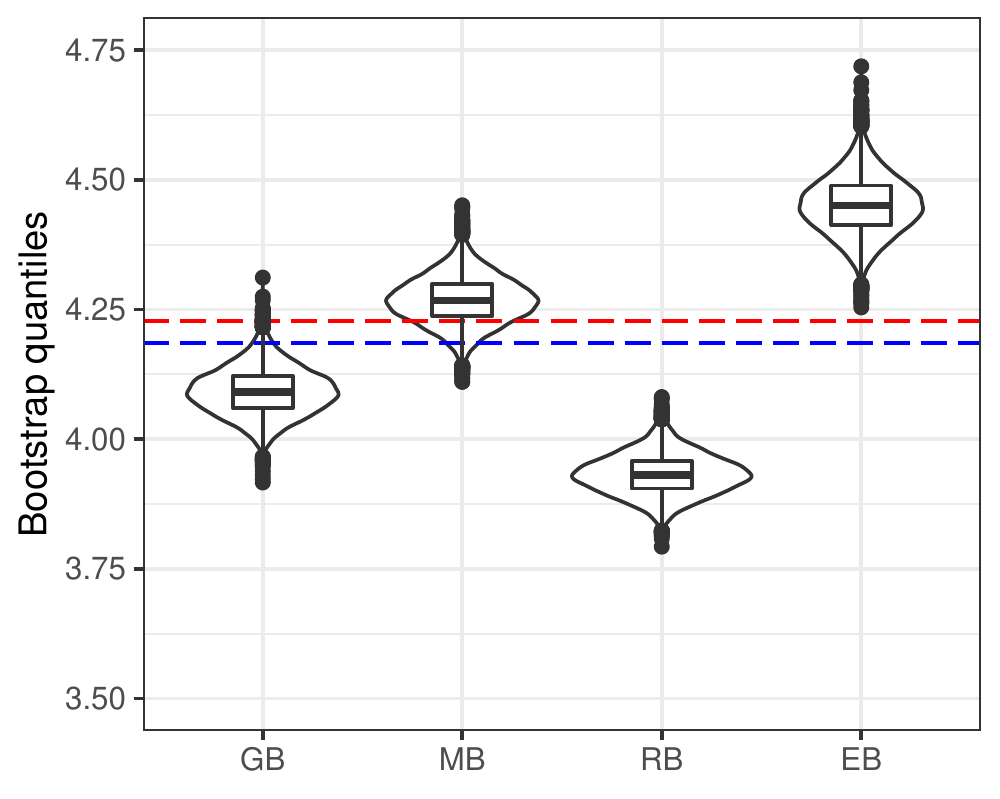}}
	\hspace{0\textwidth}
	\subfigure[$\Sigma_{j,k} = 0.2^{|j-k|}$]{
		\label{fig:exp2} 
		\includegraphics[width=0.45\textwidth]{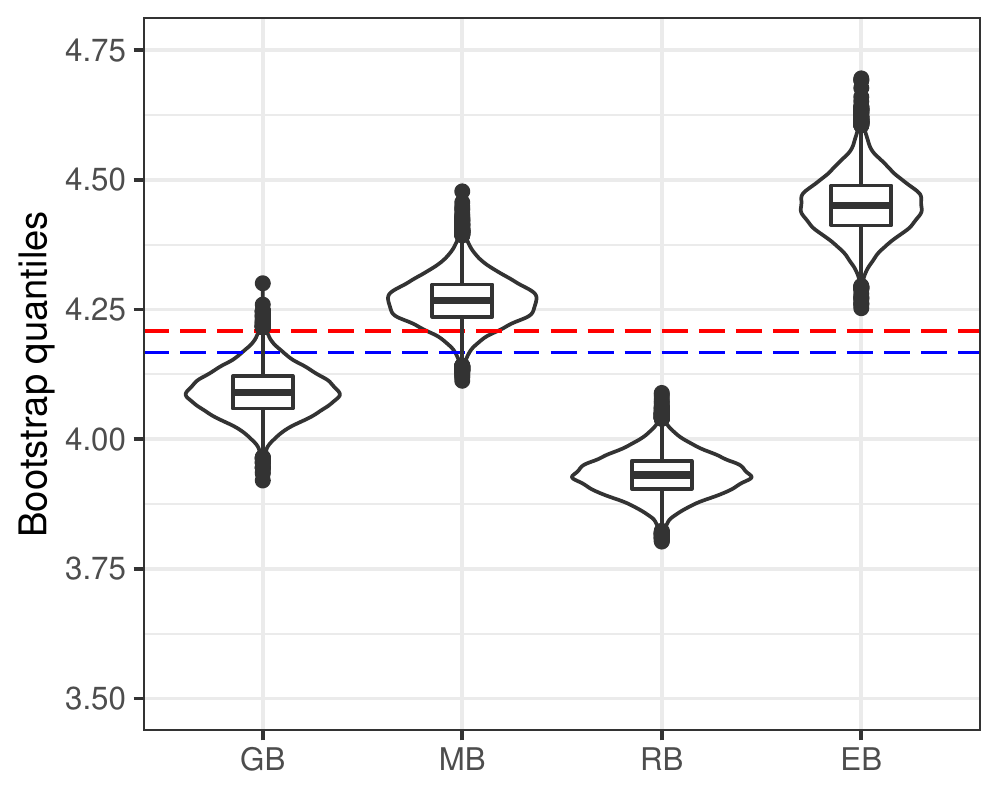}}
	\hspace{0\textwidth}
	\subfigure[$\Sigma_{j,k} = 0.8^{|j-k|}$]{
		\label{fig:exp3} 
		\includegraphics[width=0.45\textwidth]{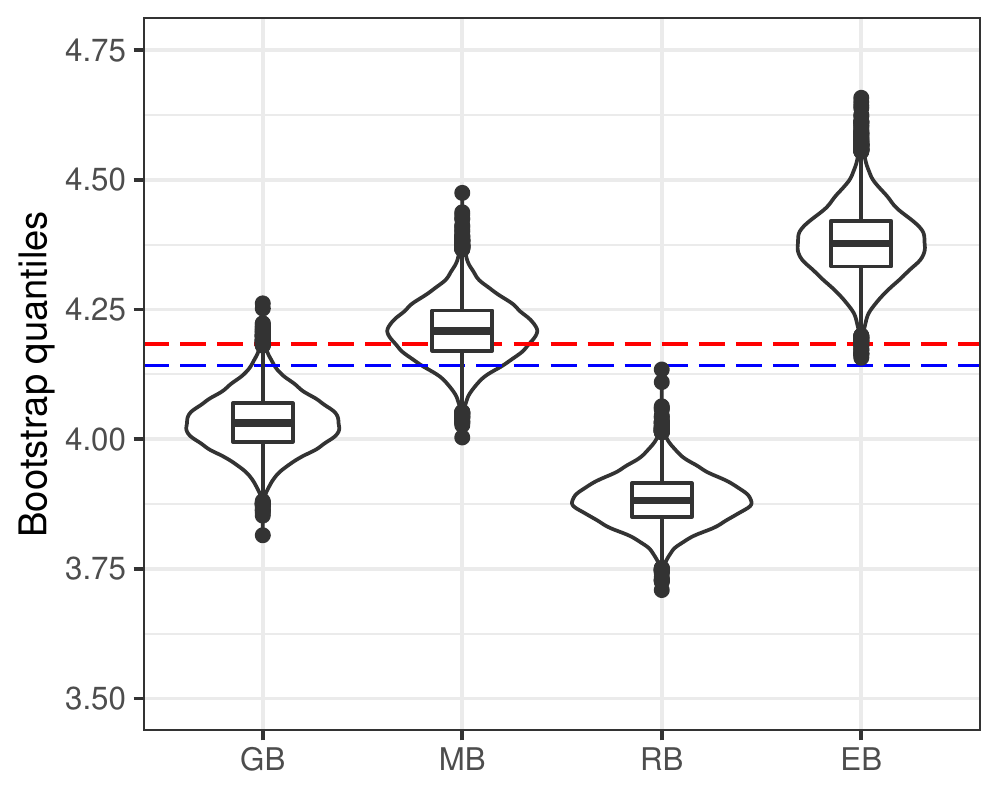}}
	\hspace{0\textwidth}
	\subfigure[$\Sigma_{j,k} = 0.2 + 0.8 I\{j=k\}$]{
		\label{fig:exp4} 
		\includegraphics[width=0.45\textwidth]{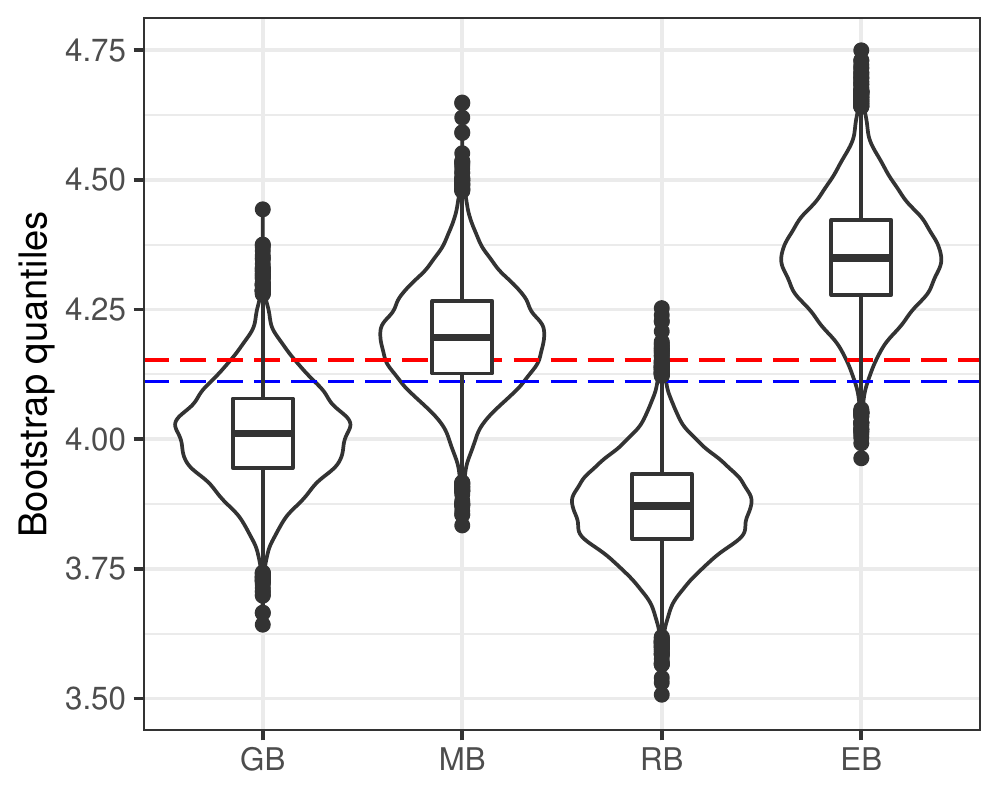}}
	\caption{Box and Violin (mirrored density) plots of the bootstrap quantiles $\{ (t_{\alpha}^*)^{(k)}, k \le K\}$ for different bootstrap schemes. (red dashed line: true $t_{\alpha}$; blue dashed line: $t_{\alpha}/1.01$)}
	\label{fig:bootstrap_quantiles}
\end{figure}

As the coverage probabilities results in Table \ref{tab:relative_frequency} do not give much information on how the conservative procedure affects the bootstrap accuracy, we may look at the quantiles $t_{\alpha}^*$ directly. The most ideal scenario of using exact bootstrap would be the simulated $(t_{\alpha}^*)^{(k)}$'s all highly concentrate on the true $t_{\alpha}$, and the conservative bootstrap then aims to make $1.01(t_{\alpha}^*)^{(k)}$'s concentrate slightly above $t_{\alpha}$, or equivalently, $(t_{\alpha}^*)^{(k)}$'s slightly above $t_{\alpha}/1.01$. The box plots and violin plots (mirrored density plots) of the bootstrap quantiles $\{(t_{\alpha}^*)^{(k)}, 1\le k \le K\}$ in Experiments (a)-(d) are given in Figure \ref{fig:bootstrap_quantiles}, where the true quantiles $t_{\alpha}$'s of the maxima $T_n$ in Experiments (a)-(d) are simulated from $5 \times 10^4$ simulations. 

In Experiments (a)-(d), it seems that the Mammen's multiplier bootstrap have the most stable performance. For this bootstrap, we observe that the slightly inflated bootstrap quantile $1.01 t_{\alpha}^*$ does make a difference in that there exist a significant amount of $(t_{\alpha}^*)^{(k)}$ in all experiments that are slightly smaller than $t_{\alpha}$ but become slightly above $t_{\alpha}$ after the small inflation, resulting a mild conservative coverage behavior. The empirical bootstrap in these settings are quite conservative by itself, but the small inflation does not really make it worse either as the coverage probabilities in Table \ref{tab:relative_frequency} only increase at most 0.3\%. The Gaussian multiplier bootstrap and the Rademacher multiplier bootstrap in general do not perform as well as the others, although the inflated quantiles do increase their coverage probabilities a little to relieve the severe under-coverage. 

Table \ref{tab:relative_frequency} demonstrates that the inflation factor $\epsilon_0 = 0.01$ translates into a roughly 1\% increase in coverage probability under moderate $p=10^3$. When the targeted coverage probability is $0.9$ or smaller, more simulation results (not given here) show that a greater 2\% or even 3\% increase may be observed. This means the slightly conservative bootstrap procedure is actually relatively significant in practice; see also Figure \ref{fig:bootstrap_quantiles} for more evidence. However, when we do observe high correlation among covariates, we may use a larger $\epsilon_0$. For example, when we consider another experiment in which $\Sigma_{j,k} = 0.8 + 0.2 I\{j = k\}$, it seems $\epsilon_0 = 0.05$ or even $\epsilon_0 = 0.1$ is more reasonable. See Table \ref{tab:extra}.

\begin{table}[htb]
\centering
	\begin{tabular}{c|c c c c}
		\hline
		& $\mathbb{P}\{ T_n \le t_{\alpha}^*\}$ & $\mathbb{P}\{ T_n \le 1.01t_{\alpha}^*\}$ & $\mathbb{P}\{ T_n \le 1.05t_{\alpha}^*\}$ & $\mathbb{P}\{ T_n \le 1.1t_{\alpha}^*\}$\\
		\hline
		MB & 0.9195 & 0.9246 & 0.9401 & 0.9550\\
		\hline
		EB & 0.9262 & 0.9303 & 0.9454 & 0.9588\\
		\hline
	\end{tabular}
	\vspace{0.5ex}
	\caption{Simulated relative frequencies under $\Sigma_{j,k} = 0.8 + 0.2 I\{j = k\}$.}
	\label{tab:extra}
\end{table}

\bibliographystyle{apalike}
\bibliography{conservative}

\appendix
\renewcommand{\thesection}{\Alph{section}}


\clearpage
\section{Proofs in Section 2}

\subsection{Proof of Theorems \ref{thm:main} and \ref{thm:main-full}} \label{proof:thm:main}

Theorem \ref{thm:main-full} follows from Lemma \ref{lm:lp-to-consistency} and the following proposition. Theorem \ref{thm:main} is a corollary of Theorem \ref{thm:main-full} with $\epsilon_n = 0.01$.

\begin{proposition}\label{prop:LP}
Under the same conditions of Theorem \ref{thm:main}, it holds with at least probability $1 - 2/(np) - q_0(M, \varepsilon)$ that
\begin{align}\label{1-prop:LP-empi}
\sup_t \eta_n^*(\varepsilon, t)  &\le \overline{\eta}_n(\varepsilon) := C_{c_0, \mathrm{boot}} \min \bigg\{\frac{(\log (np))^{3/2}}{n^{1/2}}\frac{M^2}{\varepsilon^2}, \frac{(\log (np))^{3/2}(\log p)^{1/2}}{n^{1/2}}\frac{M^2}{\varepsilon \overline{\sigma}},
\\ \nonumber
&\hspace{12em} \Big(\frac{(\log(np))^3(\log p)^2}{n} \Big)^{1/4} \frac{M}{\overline{\sigma}} \bigg\} + 2/(np) +  3q_0(M, \varepsilon). 
\end{align}
\end{proposition}

Define $\eta_1 := \overline{\eta}_n(\epsilon_n t_{\alpha + \eta_0}/(1 + \epsilon_n))$. The constant $C_{\mathrm{boot}, c_0, \eta_0}$ in \eqref{bound-thm:main} can be large enough so that a nontrivial bound yields $\eta_1 \le \eta_0$. We have
\begin{align*}
\eta_{n, \alpha}^* &\le \mathbb{P} \Big\{ \eta_n^* \Big( \frac{\epsilon_n t_{\alpha + \eta_1} }{1 + \epsilon_n} , t_{\alpha + \eta_1} \Big) \ge \eta_1 \Big\} + \eta_1
\cr
&\le \mathbb{P} \Big\{ \eta_n^* \Big( \frac{\epsilon_n t_{\alpha + \eta_0} }{1 + \epsilon_n} , t_{\alpha + \eta_1} \Big) \ge \eta_1 \Big\} + \eta_1
\cr
&\le \frac{2}{np} + q_0 \Big(M, \frac{\epsilon_n t_{\alpha + \eta_0} }{1 + \epsilon_n} \Big) + \overline{\eta}_n \Big(\frac{ \epsilon_n t_{\alpha + \eta_0}}{1 + \epsilon_n} \Big),
\end{align*}
where the first inequality follows from Lemma \ref{lm:lp-to-consistency}, the second from the fact that $\sup_t \eta_n^*(\varepsilon, t)$ is non-increasing w.r.t $\varepsilon$ and the third from Proposition \ref{prop:LP}. The proof is complete with the following proof of Proposition~\ref{prop:LP}.

\medskip
{\sc Proof of Proposition \ref{prop:LP}.} We first prove it for the empirical bootstrap in {\sc Part 1}. The multiplier bootstrap case is similar, so we only point out some differences in {\sc Part 2}. We assume $\mathbb{E} X_{i,j} = 0$ for all $i, j$.

\smallskip
{\sc Part 1: Empirical Bootstrap.} Let $u_n = M (n/(\log(np))^{1/4}$ and $a_n = \sqrt{n} \varepsilon /\log (np)$. Let $\widetilde{X}_{i,j} = X_{i,j}I\{|X_{i,j}| \le u_n \} - \mathbb{E} X_{i,j}I\{|X_{i,j}| \le u_n \}$, and $\{\widetilde{X}_{i,j}^*, 1\le i \le n \}$ be uniformly sampled from $\big\{ \widetilde{X}_{i,j}  - n^{-1} \sum_{k=1}^n \widetilde{X}_{k,j}, 1\le i \le n \big\}$. Let $\widetilde{T}_n = \max_j \sum_{i=1}^n \widetilde{X}_{i,j}/\sqrt{n}$ and $\widetilde{T}_n^* = \max_j \sum_{i=1}^n \widetilde{X}_{i,j}^*/\sqrt{n} $. It follows from the pseudo triangle inequality of the LP pre-distance in \eqref{triangle_inequality} that
\begin{align}\label{0-1-pf:LP-empi}
\eta_n^*(\varepsilon, t) &\le \sup_{t} \eta_n^*(\varepsilon/2, t; \widetilde{T}_n, \widetilde{T}_n^* ) + \mathbb{P} \Big[ \big\| \sum_{i=1}^n ( X_{i,j} - \widetilde{X}_{i, j})/\sqrt{n} \Big\|_{\infty} \ge \varepsilon/4 \Big]
\\ \nonumber
& \hspace{10.5em} + \mathbb{P}^* \Big[ \big\| \sum_{i=1}^n ( X_{i,j}^* - \widetilde{X}_{i, j}^* )/\sqrt{n} \Big\|_{\infty} \ge \varepsilon/4 \Big].
\end{align}	
We bound the above three quantities on the right-hand side when
\begin{align*}
\varepsilon \ge \varepsilon_0 = \sqrt{2\,C_{c_0, 2}} \{(\log (np))^3/n\}^{1/4} M
\end{align*}
for a sufficiently large constant $C_{c_0, 2}$ to be given later. This implies that $|\widetilde{X}_{i,j}| \le 2 u_n \le 2a_n = 2\sqrt{n} \varepsilon/\log (np)$. The case of $\varepsilon < \varepsilon_0$ will be considered in the end of the proof.

\smallskip
\noindent{$\bullet$ \it Step 1.} We bound $\displaystyle \eta_n^*(\varepsilon/2, t; \widetilde{T}_n, \widetilde{T}_n^*)$ through Theorem \ref{thm:DAT} in which we let $X_i = \widetilde{X}_i^0$ be an independent copy of $\widetilde{X}_i$, $Y_i = \widetilde{X}_i^*$, $m^*=4$ and $\mathbb{E}$ be $\mathbb{E}^* = \mathbb{E} [ \cdot | \widetilde{X}_i, 1\le i \le n]$. Note the variables $\mathfrak{W}_i$ in Theorem \eqref{thm:DAT} are taken to be fixed constant $4$. It follows from Theorem \ref{thm:DAT} that
\begin{align}\label{1-1-pf:LP-empi}
&\eta_n(\varepsilon/2, t; \widetilde{T}_n, \widetilde{T}_n^*) \le K_{n, 4}(\varepsilon/2) \min \bigg\{1, \frac{\omega_n \big( 19 \varepsilon/2; \widetilde{T}_n \big) }{ \big[ 1-K_{n, 4}(\varepsilon/2) \big]_+ } \bigg\}
\end{align}
where by Theorem \ref{thm:gAC},
\begin{align}\label{1-2-pf:LP-empi}
\displaystyle \omega_n(\varepsilon; \widetilde{T}_n)  \le C_1 \bigg( \frac{(\log (np))^3 (\log p)^{1/2}}{n} \frac{M_4^4 }{\varepsilon^3 \overline{\sigma}} + \frac{\varepsilon}{\overline{\sigma}}\sqrt{\log p} \bigg).
\end{align}
We then bound $K_{n, 4}(\varepsilon/2)$. It follows from Bennett inequality that
\begin{align*}
& \sum_{m=1}^4
 \mathbb{P} \Bigg\{\bigg\| \frac{1}{n}\sum_{k=1}^n \Big(\widetilde{X}_k^{\otimes m} - \mathbb{E} \widetilde{X}_{k}^{\otimes m} \Big) \bigg\|_{\max} > (2 u_n)^mB_{n,m} \Bigg\}
\cr 
&\le \sum_{m=1}^4 (2p^m) \exp \bigg\{ - n \Big(\frac{M_4}{2 u_n}\Big)^{4 \wedge (2m)} \rho\Big((2 u_n/M_4)^{ 4 \wedge (2m)} B_{n,m}\Big) \bigg\}
\end{align*}
where $\rho(t) = (1+t)\log(1+t) - t$ satisfies $\rho(t) \ge t^2/(2K+2)$ when $t \le K$ for any $K > 2$. Recall $\log(np) \le  c_0  n$, we let $B_{n,m} = C_{c_0, 1} \big( \log (np)/n \big)^{3/4}$ for $m=1$ and $B_{n,m} = C_{c_0, 1} \log (np)/n$ for $m=2,3$ and $4$, where the constant $C_{c_0, 1}$ is large enough so that
\begin{align*}
\begin{cases}
\displaystyle \frac{4 C_{c_0, 1}^2}{ 8 C_{c_0, 1}  c_0 ^{1/4} + 2} \log(np) \ge \log(8np^{m+1}) & m=1,
\cr
\displaystyle \frac{\rho(2^4 C_{c_0, 1})}{2^4 } \log(np) \ge \log(8np^{m+1}) & m=2, 3, \hbox{ and } 4.
\end{cases}
\end{align*}
This implies that it holds with at least probability $1 - \sum_{m=1}^4 (2 p^m) \exp( - \log(8np^{m+1})) \ge 1 - (np)^{-1}$ that $\Big\| \frac{1}{n}\sum_{k=1}^n \Big(\widetilde{X}_k^{\otimes m} - \mathbb{E} \widetilde{X}_{k}^{\otimes m} \Big) \Big\|_{\max} \le (2 u_n)^mB_{n,m} $ for $m=1, \ldots, 4$. Hence, with probability at least $1-(np)^{-1}$,
\begin{align*}
d_n^{(m)} &:= \bigg\| \frac{1}{n} \sum_{k=1}^n \mathbb{E} \widetilde{X}_k^{\otimes m} - \frac{1}{n} \sum_{k=1}^n \Big(\widetilde{X}_k - \frac{1}{n}\sum_{\ell=1}^n \widetilde{X}_{\ell} \Big)^{\otimes m} \bigg\|_{\max}
\\ \nonumber
&= \Bigg\| \frac{1}{n}\sum_{k=1}^n \Big(\widetilde{X}_k^{\otimes m} - \frac{1}{n} \sum_{\ell =1}^n \mathbb{E} \widetilde{X}_{\ell}^{\otimes m} \Big)
+
\\ \nonumber
& \qquad \qquad \sum_{m_0 =1}^m {m \choose m_0} \mathrm{Sym}
\bigg[ \Big(- \frac{1}{n}\sum_{\ell=1}^n \widetilde{X}_{\ell}\Big)^{\otimes m_0} \otimes \Big( \frac{1}{n} \sum_{k=1}^n \widetilde{X}_k^{\otimes (m-m_0)}\Big)  \bigg] \Bigg\|_{\max}
\\ \nonumber
&\le \Bigg\| \frac{1}{n}\sum_{k=1}^n \Big(\widetilde{X}_k^{\otimes m} - \mathbb{E} \widetilde{X}_{k}^{\otimes m} \Big) \Bigg\|_{\max}
+ \sum_{m_0 =1}^m {m \choose m_0} \Bigg\|\frac{1}{n}\sum_{k=1}^n \Big( \widetilde{X}_{k} - \mathbb{E} \widetilde{X}_{k} \Big) \Bigg\|_{\infty}^{m_0} \times 
\\ \nonumber
& \qquad \Bigg( \bigg\| \frac{1}{n} \sum_{k=1}^n \Big(\widetilde{X}_k^{\otimes (m-m_0)} - \mathbb{E} \widetilde{X}_k^{\otimes (m-m_0)} \Big) \bigg\|_{\max} + \bigg\| \frac{1}{n} \sum_{k=1}^n \mathbb{E} \widetilde{X}_k^{\otimes (m-m_0)} \bigg\|_{\max} \Bigg)
\\ \nonumber
&\le (2 u_n)^m B_{n,m}
+ \sum_{m_0=1}^m {m \choose m_0} \big[ (2 u_n ) B_{n,1} \big]^{m_0} \times
\\ \nonumber
& \qquad \qquad \qquad \Big( (2 u_n)^{m-m_0} B_{n, m - m_0} + \big( M_4 \wedge (2 u_n) \big)^{m-m_0} \Big)
\\ \nonumber
&\le C_{c_0, m} u_n^m \log(np) /n,
\end{align*}
where $\mathrm{Sym}(A)$ denotes the symmetrization of tensor $A$ by taking the average over all permutations of the index of its elements, and then
\begin{align*}
K_{n,4}(\varepsilon/2) &\le  C_2 \bigg( \sum_{m=2}^{4} \frac{(\log (np) )^{m-1}}{n^{m/2-1}}  \frac{d_n^{(m)} }{\varepsilon^m} + \frac{(\log (np))^3 }{ n } \frac{M_4^4}{\varepsilon^4} \bigg)
\\ \nonumber
&\le C_{c_0, 2} \frac{(\log (np))^{3/2}}{n^{1/2} } \frac{M^2}{\varepsilon^2 }.
\end{align*}
As $\varepsilon \ge \varepsilon_0$ and $K_{n,4}(\varepsilon_0/2) \le 1/2$, it follows from \eqref{1-1-pf:LP-empi} and \eqref{1-2-pf:LP-empi} that with at least probability $1-(np)^{-1}$
\begin{align}\label{1-3-pf:LP-empi}
\eta_n(\varepsilon/2, t; \widetilde{T}_n, \widetilde{T}_n^*) &\le C_{c_0, 3} \min \Big\{\frac{(\log (np))^{3/2}}{n^{1/2}}\frac{M^2}{\varepsilon^2}, \frac{(\log (np))^{3/2}(\log p)^{1/2}}{n^{1/2}}\frac{M^2}{\varepsilon \overline{\sigma}}\Big\}
\\ \nonumber
&\le C_{c_0, 4} \min \Big\{\frac{(\log (np))^{3/2}}{n^{1/2}}\frac{M^2}{\varepsilon^2}, \frac{(\log (np))^{3/2}(\log p)^{1/2}}{n^{1/2}}\frac{M^2}{\varepsilon \overline{\sigma}}, 
\\ \nonumber 
&\hspace{14em} \Big(\frac{(\log(np))^3(\log p)^2}{n} \Big)^{1/4} \frac{M}{\overline{\sigma}}\Big\}.
\end{align}

\smallskip
\noindent{$\bullet$ \it Step 2.} We bound the other two terms on the right-hand side of \eqref{0-1-pf:LP-empi} using Lemmas 3 and 7 in \cite{deng2017beyond}, which yield that
\begin{align}\label{2-1-pf:LP-empi}
& \mathbb{P} \Big[ \big\| \sum_{i=1}^n ( X_{i,j} - \widetilde{X}_{i, j})/\sqrt{n} \Big\|_{\infty} \ge \varepsilon/4 \Big] \le (np)^{-1} + \mathbb{P} \big\{\| \boldsymbol{X} \|_{\max} \ge a_n \big\},
\\ \nonumber
& \mathbb{P} \Big\{ \mathbb{P}^* \Big[ \big\| \sum_{i=1}^n ( X_{i,j}^* - \widetilde{X}_{i, j}^* )/\sqrt{n} \Big\|_{\infty} \ge \varepsilon/4 \Big] > (np)^{-1} \Big\} \le (np)^{-1} + \mathbb{P} \big\{\| \boldsymbol{X} \|_{\max} \ge a_n \big\}.
\end{align}
Note that these Lemmas are applicable because $|\widetilde{X}_{i,j}| \le 2 u_n \le 2a_n$, which is essentially due to the condition $\varepsilon \ge \varepsilon_0$.

\smallskip
\noindent{$\bullet$ \it Step 3.} It follows from \eqref{0-1-pf:LP-empi}, \eqref{1-3-pf:LP-empi} and \eqref{2-1-pf:LP-empi} that for $\varepsilon \ge \varepsilon_0$,
\begin{align}\label{3-1-pf:LP-empi}
\eta_n^*(\varepsilon, t) &\le 2/(np) +  \mathbb{P} \big\{\| \boldsymbol{X} \|_{\max} \ge a_n \} + C_{c_0, 4} \min \bigg\{\frac{(\log (np))^{3/2}}{n^{1/2}}\frac{M^2}{\varepsilon^2}, 
\\ \nonumber 
&\hspace{4em} \frac{(\log (np))^{3/2}(\log p)^{1/2}}{n^{1/2}}\frac{M^2}{\varepsilon \overline{\sigma}}, \Big(\frac{(\log(np))^3(\log p)^2}{n} \Big)^{1/4} \frac{M}{\overline{\sigma}} \bigg\}
\end{align}	
with at least probability $1- 2/(np) -  \mathbb{P} \big\{\| \boldsymbol{X} \|_{\max} \ge a_n \}$. For $\varepsilon \le \varepsilon_0$, It follows from Theorem \ref{thm:gAC} and \eqref{3-2-pf:LP-empi} that
\begin{align}\label{3-2-pf:LP-empi}
\eta_n^*(\varepsilon, t) &\le \sup_{t} \big| \mathbb{P}\{T_n \le t\} - \mathbb{P}^* \{T_n^* \le t\} \big| 
\\ \nonumber
&\le \sup_t \eta_n^*(\varepsilon_0, t) + \omega_n(\varepsilon_0; T_n)
\\ \nonumber
&\le 2/(np) + 3 \mathbb{P} \Big\{\| \boldsymbol{X} \|_{\max} \ge \sqrt{2C_{c_0, 2}} u_n \Big\} + C_{c_0, 5} \Big(\frac{(\log (np))^3 (\log p)^2}{n} \Big)^{1/4} \frac{M}{\overline{\sigma}}
\end{align}
with at least probability $\displaystyle 1 - 2/(np) -  \mathbb{P} \big\{\| \boldsymbol{X} \|_{\max} \ge \sqrt{2C_{c_0, 2}} u_n \big\}$.
The final conclusion \eqref{1-prop:LP-empi} follows from \eqref{3-1-pf:LP-empi} and \eqref{3-2-pf:LP-empi} with proper rescaling on constant coefficients. The proof is complete. $\hfill\square$

\medskip
{\sc Part 2: Multiplier bootstrap.} We only point out the differences and omit repeated arguments. As we are considering multiplier bootstrap, define $\widetilde{X}_i^* = W_i \widetilde{X}_i$ for all $i$, so that \begin{align}\label{3-pf:LP-multi}
\mathbb{P} \Big\{ \mathbb{P}^* \Big[ \big\| \sum_{i=1}^n ( X_{i,j}^* - \widetilde{X}_{i, j}^* )/\sqrt{n} \Big\|_{\infty} \ge \varepsilon/4 \Big] > (np)^{-1} \Big\} \le (np)^{-1} + \mathbb{P} \big\{\| \boldsymbol{X} \|_{\max} \ge a_n \big\}
\end{align}
is instead obtained from \cite[Lemma 8]{deng2017beyond}. 

The major difference here is the way of bounding $\sup_{t} \eta_n^*(\varepsilon/2, t; \widetilde{T}_n, \widetilde{T}_n^* )$. As $W_i \widetilde{X}_{i,j}$ may be unbounded due to the general choice of sub-Gaussian $W_i$, Theorem \ref{thm:DAT} is not directly applicable. We further truncate $W_i$ as $W_i = \widetilde{W}_i  + (W_i - \widetilde{W}_i)$ where $\widetilde{W}_i = W_i I\{|W_i| \le c\log(np)\} - \mathbb{E} W_i I\{|W_i| \le c\log(np)\}$. This step is not necessary when $W_i$ are bounded. Let $\widetilde{T}_n^{**} = \max_j \sum_{i=1}^n \widetilde{W}_i \widetilde{X}_{i,j}/\sqrt{n}$. Since 
\begin{align}\label{1-pf:LP-multi}
\sup_{t} \eta_n^*(\varepsilon/2, t; \widetilde{T}_n, \widetilde{T}_n^* ) &\le \sup_{t} \eta_n^*(\varepsilon/4, t; \widetilde{T}_n, \widetilde{T}_n^{**} ) + \mathbb{P}\bigg\{  \Big\| \frac{1}{\sqrt{n}} \sum_{i=1}^n (W_i - \widetilde{W_i}) \widetilde{X}_i \Big\|_{\infty} \ge \varepsilon/4 \bigg\}.
\end{align}	
We can now apply Theorem \ref{thm:DAT} with $\mathfrak{W}_i = 2(\widetilde{W}_i \vee 1)$ to bound $\sup_{t} \eta_n^*(\varepsilon/4, t; \widetilde{T}_n, \widetilde{T}_n^{**} )$. It yields that
\begin{align}\label{2-pf:LP-multi}
&\eta_n(\varepsilon/4, t; \widetilde{T}_n, \widetilde{T}_n^{**}) \le K_{n, 4}(\varepsilon/4) \min \bigg\{1, \frac{\omega_n \big( (16c+3) \varepsilon/4; \widetilde{T}_n \big) }{ \big[ 1-K_{n, 4}(\varepsilon/2) \big]_+ } \bigg\}
\end{align}
where the bound for $\omega_n(\varepsilon; \widetilde{T}_n)$ is still as in \eqref{1-2-pf:LP-empi}.
Note that in multiplier bootstrap with third moment match, $\mathbb{E} W_i^2 = \mathbb{E} W_i^3 = 1$, so that
\begin{align*}
K_{n,4}(\varepsilon) &\le  C_{\tau_0,1} \bigg( \sum_{m=2}^{3} \frac{(\log (np) )^{m-1}}{n^{m/2-1} \varepsilon^m } \Big\| \frac{1}{n}\sum_{i=1}^n \mathbb{E} \widetilde{X}_i^{\otimes m} - \frac{1}{n} \sum_{i=1}^n (\mathbb{E} \widetilde{W}_i^m) \widetilde{X}_i^{\otimes m} \Big\|_{\max} 
\\ \nonumber
& \hspace{6em} + \frac{(\log (np))^3 }{ n \varepsilon^4} \Big( M_4^4 + \Big\| \frac{1}{n}\sum_{i=1}^n (\mathbb{E} \widetilde{W}_i^m) \widetilde{X}_i^{\otimes m} \Big\|_{\max} \bigg)
\\ \nonumber
&\le C_{\tau_0,1} \bigg( \sum_{m=2}^{4} \frac{(\log (np) )^{m-1}}{n^{m/2-1} \varepsilon^m } \Big\| \frac{1}{n}\sum_{i=1}^n (\mathbb{E} W_i^m)\mathbb{E} \widetilde{X}_i^{\otimes m} - \frac{1}{n} \sum_{i=1}^n (\mathbb{E} \widetilde{W}_i^m) \widetilde{X}_i^{\otimes m} \Big\|_{\max} 
\\ \nonumber
& \hspace{6em} + \frac{(\log (np))^3 }{ n \varepsilon^4} \Big( M_4^4 + \Big\| \frac{1}{n}\sum_{i=1}^n (\mathbb{E} W_i^m) \mathbb{E} \widetilde{X}_i^{\otimes m} \Big\|_{\max} \bigg)
\\ \nonumber
&\le C_{\tau_0,2} \bigg( \sum_{m=2}^{4} \frac{(\log (np) )^{m-1}}{n^{m/2-1} } \frac{d_{n,1}^{(m)} + d_{n,2}^{(m)} }{\varepsilon^m} + \frac{(\log (np))^3 }{ n } \frac{M_4^4}{\varepsilon^4} \bigg),
\end{align*}
where
\begin{align*}
d_{n,1}^{(m)} &:= \Big\| \frac{1}{n} \sum_{i=1}^n (\mathbb{E} W_i^m) \widetilde{X}_i^{\otimes m} - \frac{1}{n} \sum_{i=1}^n (\mathbb{E} \widetilde{W}_i^m)\widetilde{X}_i^{\otimes m} \Big\|_{\max}, 
\\ \nonumber
d_{n,2}^{(m)} &:= \Big\| \frac{1}{n} \sum_{i=1}^n (\mathbb{E} W_i^m) \widetilde{X}_i^{\otimes m} - \frac{1}{n} \sum_{i=1}^n (\mathbb{E} W_i^m)\mathbb{E} \widetilde{X}_i^{\otimes m} \Big\|_{\max}.
\end{align*}
It follows from Cauchy-Schwarz inequality and the sub-Gaussianity of $W_i$ that 
\begin{align*}
d_{n,1}^{(m)} &= |\mathbb{E} W_1^m - \mathbb{E} \widetilde{W}_1^m| \cdot \Big\|\frac{1}{n}  \sum_{i=1}^n \widetilde{X}_i^{\otimes m} \Big\|_{\max}  
\\ \nonumber
&\le u_n^m |\mathbb{E} W_1^m - \mathbb{E} \widetilde{W}_1^m|
\\ \nonumber
&\le u_n^m \mathbb{E} \big( | W_1 - \widetilde{W}_1| \cdot | W_1^{m-1} + W_1^{m-2}\widetilde{W}_1 + \cdots + \widetilde{W}_1^{m-1}| \big)
\\ \nonumber
&\le C_{m, \tau_0, 1} u_n^m \big(\mathbb{E} | W_1 - \widetilde{W}_1|^2 \big)^{1/2} 
\\ \nonumber
&\le C_{m, \tau_0, 1} u_n^m \big(\mathbb{E} W_1^2 I\{|W_i| > c\log(np)\}\big)^{1/2} 
\\ \nonumber
&\le C_{m, \tau_0, 2} u_n^m \exp\Big(-\frac{c^2(\log(np))^2}{8\tau_0^2} \Big)
\\ \nonumber
&\le C_{m,\tau_0, 3} u_n^m /n,
\end{align*}
where the last inequality follows from a large enough $c$ that may depend on $\tau_0$. 

As $|\widetilde{X}_{k, j_1} \cdots \widetilde{X}_{k, j_m}  - \mathbb{E} \widetilde{X}_{k, j_1} \cdots \widetilde{X}_{k, j_m}  | \le 2 (2u_n)^m$ and $n^{-1}\sum_{i=1}^n \mathrm{Var} (\widetilde{X}_{k, j_1} \cdots \widetilde{X}_{k, j_m}) \le (2 u_n)^{2m-4} M_4^4$, it follows from Bernstein inequality that with at least probability $1 - (np)^{-1}$,
\begin{align*}
d_{n,2}^{(m)} &\le \sqrt{2(2 u_n)^{2m-4} M_4^4 \log(6np^{m+1})/n} + 4(2u_n)^{m} \log(6np^{m+1})/(3n) 
\\ \nonumber
&\le C_{m, 1} u_n^m \log(np) /n \quad \hbox{for } m=2, 3, 4.
\end{align*}
This eventually implies that the bound in \eqref{1-3-pf:LP-empi} also applies to $\eta_n(\varepsilon/4, t; \widetilde{T}_n, \widetilde{T}_n^*)$.

It remains to bound $\displaystyle \mathbb{P}\bigg\{  \Big\| \sum_{i=1}^n (W_i - \widetilde{W_i}) \widetilde{X}_i /\sqrt{n}\Big\|_{\infty} \ge \varepsilon/4 \bigg\}$, which can be done via
\begin{align*}
& \mathbb{P}\Big\{  \big\| \sum_{i=1}^n (W_i - \widetilde{W_i}) \widetilde{X}_i/\sqrt{n} \big\|_{\infty} \ge \varepsilon/4 \Big\}
\\ \nonumber
&\le \mathbb{P}\Big\{ \max_i \Big| W_iI\{|W_i| > c \log(np)\} - \mathbb{E} W_iI\{|W_i| > c \log(np)\} \Big| \ge \frac{\sqrt{n} \varepsilon}{4n(2u_n)} \Big\}
\\ \nonumber
&\le \mathbb{P}\{ \max_i |W_i| > c \log(np)\}
\\ \nonumber
&\le 1/(np^2)
\end{align*}
for large enough $c$, where the second inequality follows from
\begin{align*}
\frac{ \max_i \mathbb{E} W_i I\{|W_i| > c \log(np)\} }{ \sqrt{n} \varepsilon/(8nu_n)} &\le 8\sqrt{n} a_n \big(\mathbb{E} W_i^2 \big)^{1/2} \big(\mathbb{P}\{|W_i| > c\log(np) \}\big)^{1/2} / \varepsilon
\\ \nonumber
&= \big( 8n/\log(np) \big) \cdot \big(\mathbb{P}\{|W_i| > c\log(np) \}\big)^{1/2} <1.
\end{align*}
for large enough $c$.
Overall, \eqref{1-3-pf:LP-empi} still holds with at least probability $1-(np)^{-1}$ for $\sup_{t} \eta_n^*(\varepsilon/2, t; \widetilde{T}_n, \widetilde{T}_n^* )$. The rest of the proof is almost identical to the empirical bootstrap, with the only remaining difference being the extra dependence of the constants on $\tau_0$ from the sub-Gaussianity of $W_i$.

The proof of Proposition \ref{prop:LP} is complete. $\hfill\square$

\subsection{Proof of Theorem \ref{thm:main-exact}}\label{proof:thm:main-exact}

This result follows from Lemma \ref{lm:lp-to-consistency} (ii) and Proposition \ref{prop:LP} in the proof of Theorem \ref{thm:main-full} almost directly. We omit the details. $\hfill\square$

\subsection{Proof of Corollary \ref{coro:example-1}}\label{proof:coro:example-1} We specify the choice of $M$ according to the conditions. The conclusion follows from plugging the $M$ into \eqref{bound-thm:main} in Theorem \ref{thm:main} for conservative bootstrap. The proof for exact bootstrap is omitted as it is similar.

For \ref{ex-bounded}, we consider $M = \max\{M_4, \{\log(np)/n\}^{1/4} B_n\}$, so that $q_0(M, t_{\alpha + \eta_0}) = 0$. 

For \ref{ex-gaussian},  
\begin{align*}
\mathbb{P} \Big\{\| \boldsymbol{X}  - \mathbb{E}  \boldsymbol{X} \|_{\max}> M\big( n / \log (np) \big)^{1/4}\Big\} \le np \frac{\exp\big\{|X_{i,j} - \mathbb{E} X_{i,j}|^2/B_n^2 \big\}}{ \exp\big\{ (M/B_n)^2(n/\log(np))^{1/2} \big\} } \le \frac{1}{np}
\end{align*}
implies $M \le C \{(\log (np))^3/ n\}^{1/4} B_n$ for a constant $C$. Note the first inequality follows from the union bound and Markov inequality. Hence, $M = \max\{C \{(\log (np))^3/ n\}^{1/4} B_n, \, M_4 \}$. Similar calculation yields $M = \max\big\{C \{(\log (np))^5/ n\}^{1/4} B_n,\ M_4\big\}$ for \ref{ex-exp}. 

For \ref{ex-moment}, it follows from the union bound and Markov inequality that
\begin{align*}
&\mathbb{P} \Big\{\| \boldsymbol{X}  - \mathbb{E}  \boldsymbol{X} \|_{\max}> \frac{\sqrt{n} (t_{\alpha + \eta_0}/101)}{\log(np)} \Big\}
\le \frac{ \sum_{i=1}^n \mathbb{E} \max_j |X_{i,j} - \mathbb{E} X_{i,j}|^q }{ \big\{ \sqrt{n} t_{\alpha + \eta_0} \big/ (101\log(np))  \big\}^q }
\le C_q \frac{(\log (np))^q}{n^{q/2 - 1}} \frac{B_n^q}{t_{\alpha + \eta_0}^q}.
\end{align*}
Hence we can set $M = M_4$.
$\hfill\square$

\subsection{Proof of Lemma \ref{lm:lp-to-consistency}}\label{proof:lm:lp-to-consistency} 
We first prove (iii). Observe
\begin{align*}
& (1-\alpha) - \mathbb{P}\{T_n\le t^*_\alpha +  \varepsilon\}
\cr &\le (1-\alpha) - \mathbb{P}\{T_n\le t_{\alpha + \eta}, t^*_\alpha +  \varepsilon \ge t_{\alpha + \eta}\}
\cr &= (1-\alpha) - \mathbb{P}\{T_n\le t_{\alpha + \eta}\} + \mathbb{P}\{T_n \le t_{\alpha + \eta},t^*_\alpha +  \varepsilon < t_{\alpha + \eta}\} 
\cr &\le \eta + \mathbb{P}\{t^*_\alpha +  \varepsilon < t_{\alpha + \eta}\}. 
\end{align*}
By the definition of $\eta_n^*\big( \varepsilon, t)$,
\begin{align*}
& \mathbb{P}\big\{t^*_\alpha +  \varepsilon < t_{\alpha + \eta}\big\} 
\cr &\le \mathbb{P}\Big[ \mathbb{P}^*\big\{T_n^* < t_{\alpha + \eta} -  \varepsilon\big\} \ge 1 - \alpha \ge \mathbb{P}\{T_n < t_{\alpha + \eta}\} + \eta \Big] 
\cr
&= \mathbb{P}\Big[ \mathbb{P}^*\big\{T_n^* < t_{\alpha + \eta} -  \varepsilon\big\} - \mathbb{P}\{T_n < t_{\alpha + \eta}\} \ge \eta \Big] 
\cr
&\le \mathbb{P} \big\{ \eta_n^* \big( \varepsilon, t_{\alpha + \eta} ) \ge \eta \big\}
\end{align*}
The other direction follows.

Next, we prove (i). By the same proof as above, we observe
\begin{align*}
& (1-\alpha) - \mathbb{P}\{T_n\le (1+  \epsilon ) t^*_\alpha \}
\cr &\le (1-\alpha) - \mathbb{P}\{T_n\le t_{\alpha + \eta}, (1+  \epsilon )t^*_\alpha \ge t_{\alpha + \eta}\}
\cr &= (1-\alpha) - \mathbb{P}\{T_n\le t_{\alpha + \eta}\} + \mathbb{P}\{T_n\le t_{\alpha + \eta}, (1+  \epsilon)t^*_\alpha < t_{\alpha + \eta}\} 
\cr &\le \eta + \mathbb{P}\{(1+  \epsilon) t^*_\alpha  < t_{\alpha + \eta}\}. 
\end{align*}
By the definition of $\eta_n^*\big( \varepsilon, t;X, X^*\big)$,
\begin{align*}
& \mathbb{P}\big\{(1+  \epsilon) t^*_\alpha  < t_{\alpha + \eta}\big\} 
\cr &\le \mathbb{P}\Big[ \mathbb{P}^*\big\{T_n^* < t_{\alpha + \eta}/(1 +  \epsilon)\big\} \ge 1 - \alpha \ge \mathbb{P}\{T_n < t_{\alpha + \eta}\} + \eta \Big] 
\cr
&= \mathbb{P}\Big[ \mathbb{P}^*\big\{T_n^* < t_{\alpha + \eta} /(1+  \epsilon) \big\} - \mathbb{P}\{T_n < t_{\alpha + \eta}\} \ge \eta \Big] 
\cr
&\le \mathbb{P} \big\{ \eta_n^* \big( \epsilon t_{\alpha + \eta} / (1+  \epsilon), t_{\alpha+ \eta}\big) \ge \eta \big\}.
\end{align*}
As (ii) follows from (i), the proof is complete.
$\hfill\square$

\section{Proofs in Section 3}
\subsection{Proof of Lemma \ref{lm:PI}}\label{proof:lm:invariance}
Let $\mathscr{A}_{i,k} =\{(A,B): A\cup B = \{1, \ldots, n\}\setminus \{k\}, |A|=i-1, |B| = n-i\}$ and
$\mathscr{A}_i =\{(A,B): A\cup B = \{1, \ldots, n\}, |A|=i, |B| = n-i\}$.
Let $X_A=\{X_k, k\in A\}$ and $Y_B=\{Y_k, k\in B\}$. We have
$\sum_{\sigma, \sigma_i =k} f(\boldsymbol{U}_{\sigma,i},\zeta_{i,k})
= \sum_{(A,B)\in \mathscr{A}_{i,k}} c_{n,i} f(X_A,Y_B,\zeta_{i,k})$ where $c_{n,i} ={\#} \big\{\sigma: \sigma_\ell\in A\ \ \forall\ \ell<i,\sigma_i=k\big\} = (i-1)!(n-i)!$.
We observe that
\begin{align*}
\sum_{(A,B)\in \mathscr{A}_i}f(X_A,Y_B) =&\ \sum_{(A,B)\in \mathscr{A}_i,k\in A}  f(X_A,Y_B)
+ \sum_{(A,B)\in \mathscr{A}_i,k\in B} f(X_A,Y_B)I\{i<n\}
\cr =&\ \sum_{(A,B)\in \mathscr{A}_{i,k}}  f(X_A,Y_B,X_k)
+ \sum_{(A,B)\in \mathscr{A}_{i+1,k}} f(X_A,Y_B,Y_k)I\{i<n\},
\end{align*}
and therefore
\begin{align*}
& n(n!) \mathbb{A}_{\sigma, i} I\{\sigma_i = k\} q_{n,i} f(\boldsymbol{U}_{\sigma,i}, \zeta_{i,k})
\cr
& = \mathbb{E} \bigg[\sum_{i=1}^n  q_{n,i}  \sum_{\sigma,\sigma_i=k} f(\boldsymbol{U}_{\sigma,i},\zeta_{i,k})\bigg|X_k, Y_k, k\le n\bigg]
\cr
&= {\sum_{i=1}^n c_{n,i}  q_{n,i}  \sum_{(A,B)\in \mathscr{A}_{i, k}} \mathbb{E} \Big[f(X_A,Y_B,\zeta_{i,k})\Big|X_k, X_k^*, k\le n\Big]}
\cr 
&= \sum_{i=1}^n c_{n,i} q_{n,i}  \sum_{(A,B)\in \mathscr{A}_{i, k}}\Big\{\theta_{n,i}  f(X_A,Y_B,X_k)
+ (1- \theta_{n,i} ) f(X_A,Y_B,Y_k)\Big\}
\cr 
&= \sum_{i=1}^n c_{n,i}  q_{n,i}  \theta_{n,i} \sum_{(A,B)\in \mathscr{A}_{i,k}}  f(X_A,Y_B,X_k)
\cr
& \qquad \qquad \qquad + \sum_{i=0}^{n-1} c_{n,i+1} q_{n, i+1} (1- \theta_{n, i+1}) \sum_{(A,B)\in \mathscr{A}_{i+1,k}} f(X_A,Y_B,Y_i)
\cr
&= \sum_{i=1}^n c_{n,i} q_{n,i} \theta_{n,i} \sum_{(A,B)\in \mathscr{A}_i}f(X_A,Y_B) + c_{n,1}q_{n,1}(1- \theta_{n,1}) f(Y_1,\ldots, Y_n).
\end{align*}
The last equality follows from $(n-i)q_{n,i} \theta_{n,i} = i q_{n,i+1}(1- \theta_{n,i+1})$ $\forall \ 1\le i \le n-1$ which implies $c_{n,i}  q_{n,i}  \theta_{n,i} = c_{n,i+1} q_{n, i+1} (1- \theta_{n, i+1})$ for $1\le i \le n-1$.
The proof of Lemma \ref{lm:PI} is complete. $\hfill\square$

\subsection{Proof of Theorems \ref{thm:weak_comparison} and \ref{thm:strong_comparison}}\label{proof:thm:weak_comparison} \label{proof:thm:strong_comparison}
 It follows from Taylor series expansion that $\mathrm{Rem}_1 = I + II$ in Theorem \ref{thm:weak_comparison} and $\mathrm{Rem}_2 = I + II + III$ in Theorem \ref{thm:strong_comparison}, where
\begin{align*}
I &= \sum_{i=1}^n q_{n,i} \mathbb{A}_{\sigma} \int_{0}^1 \mathbb{E}  \Big\langle f^{(m^*)} ( \boldsymbol{U}_{\sigma, i}, \tau X_{\sigma_i}), \frac{(1-\tau)^{m^*-1}}{(m^*-1)!}X_{\sigma_i}^{\otimes m^*} \Big\rangle d\tau,
\cr
II &= \sum_{i=1}^n q_{n,i} \mathbb{A}_{\sigma} \int_{0}^1 \mathbb{E}  \Big\langle f^{(m^*)} ( \boldsymbol{U}_{\sigma, i}, \tau Y_{\sigma_i}), \frac{(1-\tau)^{m^*-1}}{(m^*-1)!}Y_{\sigma_i}^{\otimes m^*}\Big\rangle d\tau,
\cr
III &= \sum_{i=1}^n q_{n,i} \mathbb{A}_{\sigma} \sum_{m=2}^{m^*-1} \frac{1}{m!}\Big\langle \mathbb{E} f^{(m)} ( \boldsymbol{U}_{\sigma, i}, 0) - \mathbb{E} f^{(m)} ( \boldsymbol{U}_{\sigma, i}, \zeta_{i, \sigma_i}), \mathbb{E} X_{\sigma_i}^{\otimes m} - \mathbb{E} Y_{\sigma_i}^{\otimes m} \Big\rangle.
\end{align*}
$\mathbb{A}_{\sigma, k}$ the operator that takes average over $k$ and $\sigma$ and expectation with respect to $\zeta_{k,i}$ only, that is,
\begin{align*}
\mathbb{A}_{\sigma, k} \big\{ h(( \boldsymbol{X} , \boldsymbol{Y})_{\sigma,k}, \zeta_{k,i}) \big\} = \textstyle n^{-1} \sum_{k=1}^n \mathbb{A}_{\sigma} \mathbb{E} \big[ h(( \boldsymbol{X} , \boldsymbol{Y})_{\sigma,k}, \zeta_{k,i}) \big|  \boldsymbol{X} , \boldsymbol{Y}, \sigma, k\big].
\end{align*}

We first bound $I$. It follows from the Stability Condition \ref{cond:SC} that
\begin{align*}
|I| 
&\le \sum_{i=1}^n q_{n,i} \mathbb{A}_{\sigma} \int_{0}^1 \mathbb{E}  \Big\langle |f^{(m^*)} ( \boldsymbol{U}_{\sigma, i}, \tau X_{\sigma_i})|, \frac{(1-\tau)^{m^*-1}}{(m^*-1)!}|X_{\sigma_i}|^{\otimes m^*} \Big\rangle d\tau
\cr
&\le \frac{1}{(m^*)!} \sum_{i=1}^n q_{n,i} \mathbb{A}_{\sigma} \Big\langle \mathbb{E} \overline{f}^{(m^*)} ( \boldsymbol{U}_{\sigma, i}, 0), \mathbb{E} |X_{\sigma_i}|^{\otimes m^*} g(\|X_{\sigma_i}\|) \Big\rangle
\cr
&\le \frac{1}{(m^*)!} \sum_{i=1}^n q_{n,i} \mathbb{A}_{\sigma} \Big\langle \mathbb{E} \overline{f}^{(m^*)} ( \boldsymbol{U}_{\sigma, i}, 0)/g(\|\zeta_{i,\sigma_i}\|), \frac{\mathbb{E} |X_{\sigma_i}|^{\otimes m^*} g(\|X_{\sigma_i}\|)}{\mathbb{E} [ 1/g(\|\zeta_{i,\sigma_i}\|) ]} \Big\rangle
\cr
&\le \frac{1}{(m^*)!} \sum_{i=1}^n q_{n,i} \mathbb{A}_{\sigma} \Big\langle \mathbb{E}  \overline{f}^{(m^*)}_{\max} ( \boldsymbol{U}_{\sigma, i}, \zeta_{i, \sigma_i}), \frac{\mathbb{E} |X_{\sigma_i}|^{\otimes m^*} g(\|X_{\sigma_i}\|)}{\mathbb{E} [1/g(\|X_{\sigma_i}\|)] \wedge \mathbb{E} [1/g(\|Y_{\sigma_i}\|)]} \Big\rangle
\cr
&= \frac{1}{(m^*)!} \bigg\langle \sum_{i=1}^n q_{n,i} \mathbb{A}_{\sigma} \, \mathbb{E}  \overline{f}^{(m^*)}_{\max} ( \boldsymbol{U}_{\sigma, i}, \zeta_{i, \sigma_i}), \frac{1}{n} \sum_{k=1}^n \frac{\mathbb{E} |X_k|^{\otimes m^*} g(\|X_k\|)}{\mathbb{E} [1/g(\|X_k\|)] \wedge \mathbb{E} [1/g(\|Y_k\|)]} \bigg\rangle,
\end{align*}
where the last equality follows from Lemma \ref{lm:PI}. As $II$ can be similarly bounded, the proof of Theorem \ref{thm:weak_comparison} is complete.

To prove Theorem \ref{thm:strong_comparison}, we bound $|III|$ as
\begin{align*}
|III| 
&= \bigg| \sum_{i=1}^n q_{n,i} \mathbb{A}_{\sigma} \sum_{m=2}^{m^*-1} \frac{1}{m!} \int_0^1 \Big\langle \mathbb{E} f^{(m^*)} ( \boldsymbol{U}_{\sigma, i}, t \zeta_{i, \sigma_i}),
\cr
& \qquad \qquad \qquad  \frac{(1- \tau)^{m^*-m-1}}{(m^*-m-1)!} \zeta_{i, \sigma_i}^{\otimes (m^*-m)} \otimes \Big( \mathbb{E} X_{\sigma_i}^{\otimes m} - \mathbb{E} Y_{\sigma_i}^{\otimes m} \Big) \Big\rangle d \tau\bigg|
\cr
&\le \sum_{i=1}^n q_{n,i} \mathbb{A}_{\sigma} \sum_{m=2}^{m^*-1} \frac{1}{m!} \int_0^1 \mathbb{E} \Big\langle \big| f^{(m^*)} ( \boldsymbol{U}_{\sigma, i}, t \zeta_{i, \sigma_i})\big|,
\cr
& \qquad  \frac{(1- \tau)^{m^*-m-1}}{(m^*-m-1)!} |\zeta_{i, \sigma_i}|^{\otimes (m^*-m)} \otimes \Big( \mathbb{E} |X_{\sigma_i}|^{\otimes m} + \mathbb{E} |Y_{\sigma_i}|^{\otimes m} \Big) \Big\rangle d \tau
\cr
&\le \sum_{i=1}^n q_{n,i} \mathbb{A}_{\sigma} \sum_{m=2}^{m^*-1} \frac{1}{m!(m^*-m)!} \Big\langle \mathbb{E}  \overline{f}^{(m^*)} ( \boldsymbol{U}_{\sigma, i}, 0),
\cr
& \qquad \qquad \qquad  \mathbb{E} |\zeta_{i, \sigma_i}|^{\otimes (m^*-m)}g(\|\zeta_{i, \sigma_i}\|) \otimes \Big( \mathbb{E} |X_{\sigma_i}|^{\otimes m} + \mathbb{E} |Y_{\sigma_i}|^{\otimes m} \Big) \Big\rangle
\cr
&\le \sum_{i=1}^n q_{n,i} \mathbb{A}_{\sigma} \sum_{m=2}^{m^*-1} \frac{1}{m!(m^*-m)!} \Big\langle \mathbb{E}  \overline{f}^{(m^*)}_{\max} ( \boldsymbol{U}_{\sigma, i}, \zeta_{i, \sigma_i}),
\cr
& \qquad \qquad \qquad  \frac{\mathbb{E} |\zeta_{i, \sigma_i}|^{\otimes (m^*-m)}g(\|\zeta_{i, \sigma_i}\|)}{\mathbb{E} [1/ g(\|\zeta_{i, \sigma_i}\|)]} \otimes \Big( \mathbb{E} |X_{\sigma_i}|^{\otimes m} + \mathbb{E} |Y_{\sigma_i}|^{\otimes m} \Big) \Big\rangle
\cr
&= \sum_{m=2}^{m^*-1} \frac{1}{m!(m^*-m)!} \Big\langle \sum_{i=1}^n q_{n,i} \mathbb{A}_{\sigma} \mathbb{E}  \overline{f}^{(m^*)}_{\max} ( \boldsymbol{U}_{\sigma, i}, \zeta_{i, \sigma_i}), \mu_g^{(m^*-m, m)} \Big\rangle,
\end{align*}
where by H\"older's inequality $\mu_g^{(m^*-m, m)} \le \mu_g^{(m^*)}$. Together with the bound for $|I| + |II|$ in Theorem \ref{thm:weak_comparison}, we have
\begin{align*}
|\mathrm{Rem}_2| &\le \frac{1}{(m^*)!} \Bigg\langle \sum_{i=1}^n q_{n,i} \mathbb{A}_{\sigma} \mathbb{E}  \overline{f}^{(m^*)} ( \boldsymbol{U}_{\sigma, i}, \zeta_{i, \sigma_i}), \sum_{0 \le m \le m^*-1,\, m \neq 1} { m^* \choose m } \mu_g^{(m^*-m, m)} \Bigg\rangle
\cr
& \le \frac{2^{m^*}-m^*-1}{(m^*)!} \Bigg\langle \sum_{i=1}^n q_{n,i} \mathbb{A}_{\sigma} \mathbb{E}  \overline{f}^{(m^*)}_{\max} ( \boldsymbol{U}_{\sigma, i}, \zeta_{i, \sigma_i}), \ \mu_g^{(m^*)} \Bigg\rangle.
\end{align*}
The proof of Theorem \ref{thm:strong_comparison} is complete. $\hfill\square$

\subsection{Proof of Theorem \ref{thm:DAT}}\label{proof:thm:distributional_approximation}

Let $f = f_t$ be as in \eqref{def:f} with $\beta = 2 \log(np)/\varepsilon$ and $(\overline{f}^{(m^*)}, \overline{f}_{\max}^{m^*}, g)$ as in \eqref{f_stability_bounded}. It follows from Theorem \ref{thm:strong_comparison} ($q_{n,i}=1$, $\theta_{n,i}=i/(n+1)$) and Proposition \ref{prop:H} that 
\begin{align}\label{1-pf-thm:distributional_approximation}
& \eta_n(\varepsilon, t; T_n, T_n^Y) 
\\ \nonumber
&\le \big| \mathbb{E} f(X_1, \ldots, X_n) - \mathbb{E} f(Y_1, \ldots, Y_n) \big|
\\ \nonumber
&\le \sum_{m=2}^{m^*-1} \frac{1}{m!} \Big( \sum_{i=1}^n \mathbb{A}_{\sigma} \big\| \mathbb{E} f^{(m)} ( \boldsymbol{U}_{\sigma, i}, \zeta_{i, \sigma_i}) \big\|_1 \Big) \cdot \Big\| \frac{1}{n} \sum_{i=1}^n \mathbb{E} X_i^{\otimes m} - \frac{1}{n} \sum_{i=1}^n \mathbb{E} Y_i^{\otimes m}\Big\|_{\max} 
\\ \nonumber
& \qquad \qquad  + C_{c, m^*}^0 C_{m^*} \sum_{i=1}^n \mathbb{A}_{\sigma} \big\| \mathbb{E}  \overline{f}^{(m^*)}_{\max} ( \boldsymbol{U}_{\sigma, i}, \zeta_{i, \sigma_i}) \big\|_1 \big( M_{m^*}^{m^*} + M_{m^*, Y}^{m^*} \big)
\\ \nonumber
&\le K_{n, m^*}(\varepsilon) \eta_{n, \mathbb{A}}(t),
\end{align}
where $\eta_{n, \mathbb{A}}(t) := n^{-1} \sum_{i=1}^n \mathbb{A}_{\sigma} \mathbb{P} \{ t -  (2c+1) \varepsilon\le T_{\boldsymbol{U}_{\sigma, i}, \zeta_{i, \sigma_i}} \le t + 2c\varepsilon \}$ and $T_{\boldsymbol{U}_{\sigma, i}, \zeta_{i, \sigma_i}} $ is the counterpart of $T_n$ for $(\boldsymbol{U}_{\sigma, i}, \zeta_{i, \sigma_i})$, that is, 
\begin{align*}
\textstyle T_{\boldsymbol{U}_{\sigma, i}, \zeta_{i, \sigma_i}} = \big(\sum_{k=1}^{i-1} X_{\sigma_k} + \zeta_{i, \sigma_i} + \sum_{k=i+1}^n Y_{\sigma_k} \big) \big/\sqrt{n}.
\end{align*}
Note that $f^{(m)}(\boldsymbol{U}_{\sigma, i}, \zeta_{i,\sigma_i}) = f^{(m)}(\boldsymbol{U}_{\sigma, i}, \zeta_{i, \sigma_i}) I\{ t  - \varepsilon \le T_{\boldsymbol{U}_{\sigma, i}, \zeta_{i, \sigma_i}} \le t\} = f^{(m)}(\boldsymbol{U}_{\sigma, i}, \zeta_{i, \sigma_i}) I\{ t  - (2c + 1)\varepsilon \le T_{\boldsymbol{U}_{\sigma, i}, \zeta_{i, \sigma_i}} \le t + 2c \varepsilon\} $.

 As $\mathbb{P} \{ t -  (2c+1) \varepsilon\le T_{\boldsymbol{U}_{\sigma, i}, \zeta_{i, \sigma_i}} \le t + 2c\varepsilon \} = \mathbb{P} \{ T_{\boldsymbol{U}_{\sigma, i}, \zeta_{i, \sigma_i}} \le t + 2c\varepsilon  \} - \mathbb{P} \{ T_n^Y < t + (2c+1)\varepsilon  \} + \mathbb{P} \{ T_n^Y  < t + (2c+1)\varepsilon \} - \mathbb{P} \{ T_n^Y < t - (2c+2)\varepsilon \} + \mathbb{P} \{ T_n^Y < t - (2c+2)\varepsilon \} - \mathbb{P} \{ T_{\boldsymbol{U}_{\sigma, i}, \zeta_{i, \sigma_i}} < t - (2c+1)\varepsilon  \}$, we have
\begin{align}\label{ineq-pf-thm:distributional_approximation}
\sup_t \eta_{n, \mathbb{A}}(t) &\le \sup_t \eta_{n, \mathbb{A},1}(t) + \sup_t \eta_{n, \mathbb{A},2}(t) + \omega_n\big((4c+3)\varepsilon; T_n^Y \big), \hbox{ where }
\\ \nonumber
\eta_{n, \mathbb{A}, 1}(t) &:= \textstyle n^{-1}\sum_{i=1}^n \mathbb{A}_{\sigma} \big[ \mathbb{P} \{ T_{\boldsymbol{U}_{\sigma, i}, \zeta_{i, \sigma_i}} \le t - \varepsilon \} - \mathbb{P} \{ T_n^Y < t \} \big]
\\ \nonumber
\eta_{n, \mathbb{A}, 2}(t)  &:= \textstyle n^{-1} \sum_{i=1}^n \mathbb{A}_{\sigma} \big[ \mathbb{P} \{ T_n^Y < t - \varepsilon \} - \mathbb{P} \{ T_{\boldsymbol{U}_{\sigma, i}, \zeta_{i, \sigma_i}} < t \} \big].
\end{align}

To bound $\eta_{n, \mathbb{A}, 1}(t)$, we observe
\begin{align*}
\eta_{n, \mathbb{A}, 1}(t) 
&\le \frac{1}{n} \sum_{i=1}^n \mathbb{A}_{\sigma} \Big[\mathbb{E} f (\boldsymbol{U}_{\sigma, i}, \zeta_{i, \sigma_i}) - \mathbb{E} f(Y_1, \ldots, Y_n) \Big]
\cr
&= \frac{1}{n}\sum_{i=1}^n \mathbb{A}_{\sigma} \mathbb{E} \Big\{ \frac{i}{n+1} \sum_{k=1}^i \big[ f (\boldsymbol{U}_{\sigma,k}, X_{\sigma_k}) - f (\boldsymbol{U}_{\sigma,k}, Y_{\sigma_k}) \big]
\cr
& \qquad \qquad + \frac{n+1-i}{n+1} \sum_{k=1}^{i-1} \big[ f (\boldsymbol{U}_{\sigma,k}, X_{\sigma_k} ) - f (\boldsymbol{U}_{\sigma,k}, Y_{\sigma_k}) \big] \Big\}
\cr
&= \sum_{i=1}^n \mathbb{A}_{\sigma} \frac{n+1 - i}{n+1} \big[ \mathbb{E} f (\boldsymbol{U}_{\sigma,i}, X_{\sigma_i} ) - \mathbb{E} f_{t_1} (\boldsymbol{U}_{\sigma,i}, Y_{\sigma_i}) \big].
\end{align*}
Here we apply Theorem \ref{thm:strong_comparison} again but with $q_{n,i}=(n+1-i)/(n+1)$. Let $\zeta_{i, k}^\dagger = \delta_i^\dagger X_k + (1- \delta_i^\dagger) Y_k$ where $\delta_i^\dagger \overset{ind.}{\sim} \mathrm{Bernoulli}(\theta_{n,i} = i/(n+2))$. It follows that
\begin{align}\label{3-pf-thm:distributional_approximation}
|\eta_{n, \mathbb{A}, 1}(t)| &\le K_{n, m^*}(\varepsilon) \cdot \frac{1}{n} \sum_{i=1}^n \frac{n+1-i}{n+1} \mathbb{A}_{\sigma} \mathbb{P} \big\{ t -  (2c+1) \varepsilon \le T_{\boldsymbol{U}_{\sigma, i}, \zeta_{i, \sigma_i}^\dagger} \le t +  2c \varepsilon \big\}
\\ \nonumber
&= K_{n, m^*}(\varepsilon) \cdot \frac{1}{n} \sum_{i=1}^n\mathbb{A}_{\sigma} \Big[ \frac{n+1-i}{n+1} \frac{i}{n+2} \mathbb{P} \big\{ t -  (2c+1) \varepsilon \le T_{\boldsymbol{U}_{\sigma, i}, X_{\sigma_i}} \le t + 2c \varepsilon \big\} 
\cr
&\hspace{8em}  + \frac{n+1-i}{n+1} \frac{n+2 - i}{n+2} \mathbb{P} \big\{ t -  (2c+1) \varepsilon \le T_{\boldsymbol{U}_{\sigma, i}, Y_{\sigma_i}} \le t + 2c \varepsilon \big\} \Big]
\cr
&\le K_{n, m^*}(\varepsilon) \cdot \frac{1}{n} \sum_{i=1}^n\mathbb{A}_{\sigma} \Big[ \frac{i}{n+1} \mathbb{P} \big\{ t -  (2c+1) \varepsilon \le T_{\boldsymbol{U}_{\sigma, i}, X_{\sigma_i}} \le t + 2c \varepsilon \big\} 
\cr
&\hspace{8em} + \frac{n+1-i}{n+1} \mathbb{P} \big\{ t -  (2c+1) \varepsilon \le T_{\boldsymbol{U}_{\sigma, i}, Y_{\sigma_i}} \le t + 2c \varepsilon \big\} \Big]
\\ \nonumber
&= K_{n, m^*}(\varepsilon) \cdot \eta_{n, \mathbb{A}}(t).
\end{align}
Almost identically, we can show the bound $K_{n, m^*}(\varepsilon) \cdot \eta_{n, \mathbb{A}}(t)$ also applies to $|\eta_{n, \mathbb{A}, 2}(t) |$.
 
Now let $\overline{\rho}_n(\varepsilon) := \max \{\sup_t \eta_n(\varepsilon, t; T_n, T_n^Y), \sup_t \eta_{n, \mathbb{A}, 1}(t) + \sup_t \eta_{n, \mathbb{A}, 2}(t) \}$. It follows from \eqref{1-pf-thm:distributional_approximation}, \eqref{ineq-pf-thm:distributional_approximation} and \eqref{3-pf-thm:distributional_approximation} that
$\overline{\rho}_n(\varepsilon) \le 2 K_{n, m^*}(\varepsilon) \cdot \sup_t \eta_{n, \mathbb{A}}(t) \le 2 K_{n, m^*}(\varepsilon) \min\{1, \overline{\rho}_n(\varepsilon) + \omega_n \big( (4c+3)\varepsilon ; T_n^Y \big) \}$. The final conclusion then follows as we can by symmetry replace $T_n^Y$ with $T_n$ when bounding $\eta_{n, \mathbb{A}}(t)$. $\hfill\square$

\subsection{Proof of Theorem \ref{thm:gAC}}\label{proof:thm:general_anti-concentration}
As 
\begin{align}\label{0-1-pf:gAC}
\mathbb{P} \{ t - \varepsilon \le T_n < t \} &\le \eta_n(\varepsilon, t+\varepsilon; T_n, T_n^Y) + \eta_n(\varepsilon, t+\varepsilon; T_n, T_n^Y)
\\ \nonumber
& \hspace{10em} + \mathbb{P} \{ t - 2\varepsilon \le T_n^Y < t+\varepsilon \}
\\ \nonumber
&\le 2 \sup_{t} \eta_n(\varepsilon, t; T_n, T_n^Y) + \omega_n(3 \varepsilon; T_n^Y),
\end{align}
we derive the anti-concentration bound for a general $T_n$ by that of $T_n^Y$ with Gaussian components and their LP pre-distance. To this end, let $Y_i = W_i \widetilde{X}_i$, where $W_i$ has Gaussian component and satisfies $\mathbb{E} W_i=0$ and $\mathbb{E} W_i^2 = \mathbb{E} W_i^3=1$. We give an specific choice of such $W_i$ as follows to concretize our later calculations. Let $W_i = (2/\sqrt{7}) \delta_i Z_i + 2(1- \delta_i) W_i^0$, where $\{Z_i, \delta_i, W_i^0, i=1, \ldots, n\}$ are mutually independent, $Z_i \sim \mathcal{N}(0,1)$, $\delta_i \sim \mathrm{Bernoulli}(7/8)$ and $W_i^0$ follows $\mathbb{P} \{W_i^0 = \big( 1 \pm \sqrt{5} \big)/2\} = \big( \sqrt{5} \mp 1 \big)/(2\sqrt{5})$. It is easy to show there exists a fixed constant $c_1$ such that $\mathbb{P}\{|W_i| > t\} \le c_1 \exp(-t^2/2)$.

The proof is organized as follows. Under the conditions on $a_n$:
\begin{enumerate}[label= (Cond-\arabic*), leftmargin=1in]
\item $2 M_4^4/a_n^2 \le \sigma_{(1)}^2/2$, \label{cond-1-pf:gAC}

\item $a_n^2 \log\big(j^2 \overline{\sigma}/(\varepsilon\sqrt{\log p}) \big)/n < \sigma_j^2/227$ for all $1 \le j \le p$, \label{cond-2-pf:gAC}
\end{enumerate}
we bound $\omega_n(\varepsilon; T_n^Y)$ in Step 1, bound $\eta_n(\varepsilon, t; T_n, T_n^Y)$ in Step 2 and collect their rates to obtain \eqref{1-thm:gAC} in Step 3. 
In Step 4, we show the above conditions on $a_n$ hold.

\smallskip
\noindent{$\bullet$ \it Step 1.} Let $Z_n^\dagger = \sum_{i=1}^n W_i \widetilde{X}_i/\sqrt{n}$. Denote $\mathbb{P}^\dagger = \mathbb{P}\{\cdot | (\delta_i, W_i^0, \widetilde{X}_i) \ \forall i\}$. Under $\mathbb{P}^\dagger$, $Z_n^\dagger$ is an Gaussian vector with individual means and variance $(\widetilde{\sigma}_j^\dagger)^2 = \sum_{i=1}^n \delta_i \widetilde{X}_{i,j}^2/n$. The anti-concentration of $T_n^Y$ is bounded by
\begin{align*}
\omega_n(\varepsilon; T_n^Y) 
 = \sup_t \mathbb{E} \Big[ \mathbb{P}^\dagger \big\{ t - \varepsilon \le \max_j Z_{n, j}^\dagger < t  \big\} \Big] \le \mathbb{E} \Big[ \sup_t  \mathbb{P}^\dagger \big\{ t - \varepsilon \le \max_j Z_{n, j}^\dagger < t  \big\} \Big].
\end{align*}	
Let $\widetilde{\sigma}_j^2 = \frac{1}{n} \sum_{i=1}^n \mathbb{E} \widetilde{X}_{i,j}^2$ for all $1\le j \le p$. The mean and variance of $(\widetilde{\sigma}_j^\dagger)^2$ are
\begin{align*}
\mathbb{E} \big[ (\widetilde{\sigma}_j^\dagger)^2 \big]  = (7/8) \widetilde{\sigma}_j^2 \hbox{ and } \mathrm{Var} \big[ (\widetilde{\sigma}_j^\dagger)^2 \big] = \frac{1}{n^2} \sum_{i=1}^n \mathrm{Var}(\delta_i \widetilde{X}_{i,j}^2 )\le (7/2) a_n^2 \widetilde{\sigma}_j^2/n.
\end{align*}
Note that $|\widetilde{X}_{i,j}| \le 2a_n$. It follows from \ref{cond-1-pf:gAC} that
\begin{align*}
\sigma_j^2 - \widetilde{\sigma}_j^2 &= \frac{1}{n}\sum_{i=1}^n \mathbb{E} X_{i,j}^2I\{|X_{i,j}| > a_n\} + \frac{1}{n}\sum_{i=1}^n \Big( \mathbb{E} X_{i,j} I\{|X_{i,j}| > a_n\} \Big)^2 
\cr
&\le 2 M_4^4/a_n^2 \le \sigma_j^2/2,
\end{align*}
which implies $\sigma_j^2/2 \le \widetilde{\sigma}_j^2 \le \sigma_j^2$. Bernstein's inequality yields that with $\varepsilon_0 = \varepsilon \sqrt{\log p} / \overline{\sigma}$
\begin{align*}
\mathbb{P} \bigg\{  (\widetilde{\sigma}_j^\dagger)^2 - \frac{7}{8} \widetilde{\sigma}_j^2 < - \sqrt{2\log(j^2/\varepsilon_0) \cdot 7a_n^2 \widetilde{\sigma}_j^2 /(2n)}  - 2(2a_n)^2 \log(j^2/\varepsilon_0)/(3n) \bigg\} \le \varepsilon_0 / j^2.
\end{align*}
By \ref{cond-2-pf:gAC}, $\sqrt{2\log(j^2/\varepsilon_0) \cdot 7a_n^2 \widetilde{\sigma}_j^2/(2n)}  + 8a_n^2 \log(j^2/\varepsilon_0)/(3n) < 3\sigma_j^2/16$, which implies $(\widetilde{\sigma}_j^\dagger)^2 \ge 7\widetilde{\sigma}_i^2/8 - 3 \sigma_j^2/16 \ge \sigma_j^2/4$ with at least probability $1- \sum_{j=1}^p j^{-2}  \varepsilon_0 \ge 1 - 2 \varepsilon \sqrt{\log p}/\overline{\sigma} $. It follows from \cite[Theorem 10]{deng2017beyond} that
\begin{align}\label{1-1-pf:gAC}
\omega_n(\varepsilon; T_n^Y)
&\le \mathbb{E} \Big[ \sup_t  \mathbb{P}^\dagger \big\{ t - \varepsilon \le \max_j Z_{n, j}^\dagger < t  \big\} \Big]
\\ \nonumber
&\le \frac{\varepsilon}{\overline{\sigma}} \big( 2 + \sqrt{2\log p}) + 2 \frac{\varepsilon}{\overline{\sigma}}\sqrt{\log p}
\\ \nonumber
&\le C_1 \frac{\varepsilon}{\overline{\sigma}}\sqrt{ \log p }.
\end{align}

\smallskip
\noindent{$\bullet$ \it Step 2.} Let $\widetilde{T}_n = \max_j \sum_{i=1}^n \widetilde{X}_{i,j}/\sqrt{n}$. 
Let $\widehat{W}_i = W_i I\{|W_i| \le 2\sqrt{2} \log(np)\}$, and accordingly, $\widehat{Y}_i = \widehat{W}_i \widetilde{X}_i$ and $T_n^{\widehat{Y}}$ be the counterpart of $T_n^{Y}$ for $\{\widehat{Y}_i\}$. 
By the definition of LP pre-distance \eqref{def:general_LP}, we have a pseudo-triangle inequality similar to \eqref{triangle_inequality} as 
\begin{align}\label{2-1-pf:gAC}
& \sup_t \eta_n(\varepsilon, t; T_n, T_n^Y)
\\ \nonumber
&\le \sup_t \eta_n(\varepsilon/4, t; \widetilde{T}_n, T_n^{\widehat{Y}}) + \sup_t \eta_n(\varepsilon/4, t; T_n^{\widehat{Y}}, T_n^Y) + \sup_t \eta_n(\varepsilon/2, t; T_n, \widetilde{T}_n)
\\ \nonumber
&\le \sup_t \eta_n(\varepsilon/4, t; \widetilde{T}_n, T_n^{\widehat{Y}}) + \mathbb{P} \Big[ \big\|\sum_{i=1}^n ( W_i - \widehat{W}_i) \widetilde{X}_i /\sqrt{n} \big\|_{\infty} \ge \varepsilon/4 \Big] 
\\ \nonumber
& \hspace{11em} + \mathbb{P} \Big[ \big\|\sum_{i=1}^n ( X_i - \widetilde{X}_i )/\sqrt{n} \big\|_{\infty} \ge \varepsilon/2 \Big],
\end{align}
where the second term on the right-hand side can be bounded by
\begin{align}\label{2-2-pf:gAC}
&\mathbb{P} \Big[ \big\|\sum_{i=1}^n ( W_i - \widehat{W}_i) \widetilde{X}_i /\sqrt{n} \big\|_{\infty} \ge \varepsilon/4 \Big] 
\\ \nonumber
&= \mathbb{P} \Big[ \big\|\sum_{i=1}^n \widetilde{X}_i W_i I\{|W_i| > 2\sqrt{2} \log(np) \}  /\sqrt{n} \big\|_{\infty} \ge \varepsilon/4 \Big]
\\ \nonumber
&\le \mathbb{P}\{\max_i |W_i| > 2\sqrt{2} \log(np)\}
\\ \nonumber
&\le n \exp( - 8(\log (np))^2 / 2) 
\\ \nonumber
&\le 1/(n^3p^4).
 \end{align}

Since $|\widetilde{X}_{i,j}| \vee |\widehat{Y}_{i,j}| \le \{c_0(|\widehat{W}_i| \vee 1)\} \sqrt{n} \varepsilon/(\log (np))$, $c_0(|\widehat{W}_i| \vee 1) \le 2\sqrt{2} c_0 \log(np)$ and there exists a constant $C_{c_0, 1}$ such that
\begin{align*}
\max_i \frac{\mathbb{E} \exp \big\{ 16 c_0(|\widehat{W}_i| \vee 1) \big\}}{ \mathbb{E}\exp \big\{ - 16 c_0(|\widehat{W}_i| \vee 1) \big\} } \le C_{c_0, 1},
\end{align*}
we apply Theorem \ref{thm:DAT} with $m^*=4$ to bound $\eta_n(\varepsilon, t; \widetilde{T}_n, T_n^{\widehat{Y}})$ as
\begin{align}\label{2-3-pf:gAC}
\eta_n(\varepsilon, t; \widetilde{T}_n, T_n^{\widehat{Y}}) \le K_{n, 4}(\varepsilon) \min \bigg\{1, \frac{ \omega_n \big( (8\sqrt{2}c_0 + 3)\varepsilon; T_n^{\widehat{Y}}\big)}{[1- K_{n,4}(\varepsilon)]_+} \bigg\}.
\end{align}
As $\mathbb{E} W_i^2 = \mathbb{E} W_i^3  = 1$,
\begin{align*}
K_{n, 4}(\varepsilon) &\le \sum_{m=2}^{3} C_m \frac{(\log (np))^{m-1} }{n^{m/2-1} \varepsilon^m} \Big(\mathbb{E} W_1^m I\{|W_1| > 2\sqrt{2} \log(np)\} \Big) \cdot \Big\|  \frac{1}{n}\sum_{i=1}^n \mathbb{E} \widetilde{X}_i^{\otimes m} \Big\|_{\max} 
\\ \nonumber
& \hspace{4em} + C_{c_0, 2} \frac{(\log (np))^3}{n} \frac{M_4^4}{\varepsilon^4}
\\ \nonumber
&\le \sum_{m=2}^{3} C_m \frac{(\log (np))^{m-1} }{n^{m/2-1} \varepsilon^m} \sqrt{ \mathbb{E} W_i^{2m} \cdot \mathbb{P} \{|W_1| > 2\sqrt{2} \log (np)\}} \cdot a_n^m  
\\ \nonumber
& \hspace{4em} + C_{c_0, 2} \frac{(\log (np))^3}{n} \frac{M_4^4}{\varepsilon^4}
\\ \nonumber
&\le \sum_{m=2}^{3} C_m \frac{(\log (np))^{m-1}}{n^{m/2-1} \varepsilon^m} C_2 (np)^{-2} \Big(\frac{c_0 \sqrt{n} \varepsilon }{\log (np)} \Big)^m + C_{c_0, 2} \frac{(\log (np))^3}{n} \frac{M_4^4}{\varepsilon^4}
\\ \nonumber
&\le C_{c_0, 3} \Big\{ \frac{1}{np^2\log (np)} + \frac{(\log (np))^3}{n} \frac{M_4^4}{\varepsilon^4} \Big\}
\\ \nonumber
&\le C_{c_0, 4} \frac{(\log (np))^3}{n} \frac{M_4^4}{\varepsilon^4}.
\end{align*}
The last inequality holds for $\varepsilon \ge \varepsilon_1 := (2 C_{c_0,4})^{1/4} \{(\log (np))^3/n\}^{1/4}M_4$. 

On the other hand, we know from step 1 that
\begin{align*}
& \omega_n \big( (8\sqrt{2}c_0 + 3)\varepsilon; T_n^{\widehat{Y}}\big)
\\ \nonumber
&= \sup_t \Big[ \mathbb{P}\big\{t - (8\sqrt{2}c_0 + 3)\varepsilon \le T_n^{\widehat{Y}} < t, \max_i |W_i| \le 2\sqrt{2}\log(np) \big\} 
\\ \nonumber
& \qquad \qquad - \mathbb{P} \big\{\mathbb{P}\big\{t - (8\sqrt{2}c_0 + 3)\varepsilon \le T_n^{\widehat{Y}} < t, \max_i |W_i| > 2\sqrt{2}\log(np) \big\} \Big]
\\ \nonumber
&\le \omega_n \big( (8\sqrt{2}c_0 + 3)\varepsilon; T_n^{Y} \big) + \mathbb{P} \big\{ \max_i |W_i| > 2\sqrt{2}\log(np) \big\}
\\ \nonumber
&\le C_{c_0, 5} \frac{\varepsilon}{\overline{\sigma}} \sqrt{\log p} + 1/(n^3 p^4)
\\ \nonumber
&\le C_{c_0,6} \frac{\varepsilon}{\overline{\sigma}} \sqrt{\log p}.
\end{align*}
The last inequality holds for $\varepsilon \ge \varepsilon_1$.

We now plug the above two upper bounds into \eqref{2-3-pf:gAC}. For $\varepsilon \ge \varepsilon_1$, $1 - K_{n, 4}(\varepsilon) \ge 1/2$, so that
\begin{align*}
\eta_n(\varepsilon, t; \widetilde{T}_n, T_n^{\widehat{Y}}) \le C_{c_0, 4} \frac{(\log (np))^3}{n} \frac{M_4^4}{\varepsilon^4} \min \Big\{ 1, \ 2 \cdot C_{c_0,6} \frac{\varepsilon}{\overline{\sigma}} \sqrt{\log p} \Big\},
\end{align*}
while for any $\varepsilon$, 
\begin{align*}
\eta_n(\varepsilon, t; \widetilde{T}_n, T_n^{\widehat{Y}}) &\le \sup_t \big| \mathbb{P}\{\widetilde{T}_n \le t\} - \mathbb{P}\{T_n^{\widehat{Y}} \le t\} \big|
\\ \nonumber
&\le \eta_n(\varepsilon_1, t; \widetilde{T}_n, T_n^{\widehat{Y}}) + \omega_n(\varepsilon_1; T_n^{\widehat{Y}})
\\ \nonumber
&\le C_{c_0, 4} \frac{(\log (np))^3}{n} \frac{M_4^4}{\varepsilon_1^4} \min \Big\{ 1, \ 2 \cdot C_{c_0,6} \frac{\varepsilon_1}{\overline{\sigma}} \sqrt{\log p} \Big\} + C_{c_0,6} \frac{\varepsilon_1}{\overline{\sigma}} \sqrt{\log p}
\\ \nonumber
&\le 2 C_{c_0,6} (2 C_{c_0,4})^{1/4} \frac{(\log (np))^{3/4} (\log p)^{1/2} }{n^{1/4}} \frac{M_4}{\overline{\sigma}}.
\end{align*}
Overall, it holds for any $\varepsilon$ that
\begin{align}\label{2-4-pf:gAC}
\eta_n(\varepsilon, t; \widetilde{T}_n, T_n^{\widehat{Y}}) &\le C_{c_0, 7} \min \bigg\{ \frac{(\log (np))^3}{n} \frac{M_4^4}{\varepsilon^4}, \frac{(\log (np))^3(\log p)^{1/2}}{n} \frac{M_4^4}{\varepsilon^3 \overline{\sigma}}, 
\\ \nonumber
&\hspace{14em} \frac{(\log (np))^{3/4} (\log p)^{1/2} }{n^{1/4}} \frac{M_4}{\overline{\sigma}}\bigg\}.
\end{align}

\smallskip
\noindent{$\bullet$ \it Step 3.} We obtain \eqref{1-thm:gAC} from \eqref{0-1-pf:gAC}, \eqref{1-1-pf:gAC}, \eqref{2-1-pf:gAC}, \eqref{2-2-pf:gAC} and \eqref{2-4-pf:gAC} as
\begin{align*}
& \omega_{n}(\varepsilon; T_n)
\cr
&\le 2 \sup_t \eta_n(\varepsilon, t; T_n, T_n^Y) + 3\omega_n(\varepsilon; T_n^Y)
\cr
&\le 2C_{c_0,7} \min \bigg\{ \frac{(\log (np))^3}{n} \frac{M_4^4}{\varepsilon^4}, \frac{(\log (np))^3(\log p)^{1/2}}{n} \frac{M_4^4}{\varepsilon^3 \overline{\sigma}}, \frac{(\log (np))^{3/4} (\log p)^{1/2} }{n^{1/4}} \frac{M_4}{\overline{\sigma}}\bigg\}
\cr
& \qquad + 2\mathbb{P} \bigg[ \Big\| \frac{1}{\sqrt{n}}\sum_{i=1}^n \big( X_{i,j} - \widetilde{X}_{i, j} \big) \Big\|_{\infty} \ge \frac{\varepsilon}{2} \bigg] + \frac{2}{n^3p^4} + 3C_1 \frac{\varepsilon}{\overline{\sigma}}\sqrt{\log p}
\cr
&\le C_{c_0} \bigg[ \frac{(\log (np))^3(\log p)^{1/2}}{n} \frac{M_4^4}{\varepsilon^3 \overline{\sigma}}  + \frac{\varepsilon}{\overline{\sigma}}\sqrt{\log p} \bigg] + 2\mathbb{P} \bigg[ \Big\| \frac{1}{\sqrt{n}}\sum_{i=1}^n \big( X_{i,j} - \widetilde{X}_{i, j} \big) \Big\|_{\infty} \ge \frac{\varepsilon}{2} \bigg].
\end{align*}
We then show \eqref{2-thm:gAC}.
\begin{align*}
& \mathbb{P} \bigg[ \Big\| \frac{1}{\sqrt{n}}\sum_{i=1}^n \big( X_{i,j} - \widetilde{X}_{i, j} \big) \Big\|_{\infty} \ge \frac{\varepsilon}{2} \bigg]
\cr
&= \mathbb{P} \bigg[ \Big\| \frac{1}{\sqrt{n}}\sum_{i=1}^n \big( X_{i,j}I\{|X_{i,j}| > a_n\} - \mathbb{E} X_{i,j}I\{|X_{i,j} > a_n\} \big) \Big\|_{\infty} \ge \frac{\varepsilon}{2} \bigg]
\cr
&\le \mathbb{P} \bigg[ \Big\| \frac{1}{\sqrt{n}}\sum_{i=1}^n X_{i,j}I\{|X_{i,j}| > a_n\} \Big\|_{\infty} \ge \frac{\varepsilon}{4} \bigg]
\cr
&\le \frac{4\sqrt{n}}{\varepsilon} a_n^{-3} \, \mathbb{E} \max_j \frac{1}{n}\sum_{i=1}^n X_{i,j}^4I\{|X_{i,j}| \ge a_n\}
\cr
&\le \frac{4(\log (np))^3}{c_0^3 n \varepsilon^4 } \mathbb{E} \max_j \frac{1}{n}\sum_{i=1}^n X_{i,j}^4I\{|X_{i,j}| \ge a_n\},
\end{align*}
where the first inequality is due to
\begin{align*}
& \max_j \Big| \frac{1}{\sqrt{n}}\sum_{i=1}^n \mathbb{E} X_{i,j}I\{ |X_{i,j}| > a_n\} \Big| \le \sqrt{n} a_n^{-3} \, \max_j \frac{1}{n}\sum_{i=1}^n \mathbb{E} X_{i,j}^4I\{|X_{i,j}| \ge a_n\} \le \varepsilon/4.
\end{align*}   

\smallskip
\noindent{$\bullet$ \it Step 4.} In this final step, we show the conditions of $a_n$, \ref{cond-1-pf:gAC} and \ref{cond-2-pf:gAC} hold as otherwise the bound in \eqref{1-thm:gAC} (also in Step 3) is trivial. Recall that $\overline{\sigma}/\sqrt{\log p} \le \sigma_{(1)} \le \overline{\sigma}$. To show \ref{cond-1-pf:gAC}, we shall have
\begin{align*}
\frac{2 M_4^4/a_n^2 }{\sigma_{(1)}^2/2} 
&\le \frac{4(\log (np))^2 (\log p)}{c_0^2 n \varepsilon^2 \overline{\sigma}^2} M_4^4 
\cr
& = \frac{4}{\log p} \Big( \frac{(\log (np))^2 (\log p)^{3/2}}{n} \frac{M_4^4 }{\varepsilon^3 \overline{\sigma}}\Big) \Big( \frac{\varepsilon}{\overline{\sigma}}\sqrt{\log p} \Big) \le 1
\end{align*}
for a sufficiently large $C_{c_0}$ as otherwise the bound \eqref{1-thm:gAC} is trivial. Similarly, to show \ref{cond-2-pf:gAC}, we look at
\begin{align*}
& \frac{a_n^2 \log\big(j^2 \overline{\sigma}/(\varepsilon\sqrt{\log p}) \big)/n}{\sigma_j^2/227} 
\cr
&\le 227 c_0^2 \frac{\varepsilon^2 \log\big(p^2 \overline{\sigma}/(\varepsilon\sqrt{\log p}) \big)}{\overline{\sigma}^2 (\log (np))^2} (\log p) 
\cr
&\le \frac{227 c_0^2}{ (\log (np)
)^2} \Big[ 2 \big(\frac{\varepsilon}{\overline{\sigma}}\sqrt{\log p}\big)^2 + \frac{\log\big(\overline{\sigma}/(\varepsilon \sqrt{\log p})\big)}{[\overline{\sigma}/(\varepsilon \sqrt{\log p})]^2} \Big]
\cr
&< 1.
\end{align*}

The proof is complete. $\hfill\square$

\end{document}